\documentclass[12pt]{amsart}
\usepackage{amssymb}
\usepackage{pstcol,pst-3d,pst-grad}

\let\oldlabel=\label
\def\prellabel{\marginparsep=1em
    \def\label##1{\oldlabel{##1}\ifmmode\else\ifinner\else
         \marginpar{{\footnotesize\ \\ \tt
                    ##1}}\fi\fi}}

\input diagrams
\diagramstyle[amstex]
\newarrow{Equal} =====
\newarrow{Into} C--->
\newarrow{Onto} ----{>>}

\def\Ker{{\operatorname{Ker}}}

\def\Hom{{\operatorname{Hom}}}

\def\ggr.aut{{\operatorname{gr.aut}}}

\def\Gr.Hom{{\operatorname{Gr.Hom}}}
\def\L{\operatorname{L}}

\def\int{{\operatorname{int}}}
\def\conv{\operatorname{conv}}
\def\het{\operatorname{ht}}
\def\E{\operatorname{E}}

\def\St{\operatorname{{\Bbb S}t}}

\def\path{\operatorname{path}}

\def\Aut{\operatorname{Aut}}
\def\End{\operatorname{End}}

\def\gp{\operatorname{gp}}

\def\Ker{{\operatorname{Ker}}}

\def\Proj{\operatorname{Proj}}

\def\Im{\operatorname{Im}}
\def\Proj{\operatorname{Proj}}
\def\Spec{\operatorname{Spec}}
\def\Hom{\operatorname{Hom}}
\def\Vect{\operatorname{Vect}}
\def\Pol{\operatorname{Pol}}
\def\GL{\operatorname{GL}}

\def\Col{\operatorname{Col}}
\def\Aff{\operatorname{Aff}}
\def\Coker{\operatorname{Coker}}

\def\CB{\mathcal{CB}}

\def\Q{\qedsymbol\kern1pt}
\def\sqq#1#2{{\hbox{\kern 0.5pt\vbox{\hbox{\kern#1ex\vrule width #2pt height #1ex}%
     \hrule height #2pt}\kern 1.5pt}}}
\def\sq{{\mathchoice{\sqq1{0.5}}{\sqq1{0.5}}{\sqq{0.9}{0.4}}
                    {\sqq{0.7}{0.3}}}}
\def\sqd{\sq}
\def\vt{{\kern1pt|}}
\def\RR{{\mathbb R}}
\def\CC{{\mathbb C}}

\def\ZZ{{\mathbb Z}}
\def\NN{{\mathbb N}}
\def\TT{{\mathbb T}}

\def\EE{{\mathbb E}}
\def\AA{{\mathbb A}}

\def\FF{{\mathbb F}}

\def\Pp{{\mathcal P}}

\def\Pp{{\mathfrak P}}

\def\Qq{{\mathfrak Q}}
\def\ee{{\mathbf e}}
\let\la=\langle
\let\ra=\rangle

\let\epsilon=\varepsilon
\let\phi=\varphi
\let\theta=\vartheta

\let\cal\mathcal
\let\Bbb\mathbb
\let\frak=\mathfrak

\advance\headheight1.15pt

\psset{unit=0.7cm,linewidth=0.4pt,arrowsize=1.5pt 4}
\def\vertex{\pscircle[fillstyle=solid,fillcolor=black]{0.0714}}
\newpsstyle{fatline}{linewidth=1pt}
\newpsstyle{fyp}{fillstyle=solid,fillcolor=verylight}
\definecolor{verylight}{gray}{0.95}
\definecolor{light}{gray}{0.9}
\definecolor{medium}{gray}{0.85}

\newtheorem{lemma}{Lemma}[section]
\newtheorem{corollary}[lemma]{Corollary}
\newtheorem{theorem}[lemma]{Theorem}
\newtheorem{proposition}[lemma]{Proposition}

\theoremstyle{definition}
\newtheorem{definition}[lemma]{Definition}
\newtheorem{remark}[lemma]{Remark}
\newtheorem{example}[lemma]{Example}

\begin{document}

\title[Polyhedral $K_2$]{Polyhedral $K_2$}

\author{Winfried Bruns \and Joseph Gubeladze}
\address{Universit\"at Osnabr\"uck,
FB Mathematik/Informatik, 49069 Osnabr\"uck, Germany}
\email{Winfried.Bruns@mathematik.uni-osnabrueck.de}
\address{A. Razmadze Mathematical Institute, Alexidze St. 1, 380093
Tbilisi, Georgia}
\email{gubel@rmi.acnet.ge}

\thanks{The second author was supported by the Deutsche
Forschungsgemeinschaft, INTAS grant
99-00817 and TMR grant ERB FMRX CT-97-0107}

\subjclass{14L27, 14M25, 19C09, 52B20}

\begin{abstract}
Using elementary graded automorphisms of polytopal algebras
(essentially the coordinate rings of projective toric varieties)
polyhedral versions of the group of elementary matrices and the
Steinberg and Milnor groups are defined. They coincide with the
usual $K$-theoretic groups in the special case when the polytope
is a unit simplex and can be thought of as compact/polytopal
substitutes for the tame automorphism groups of polynomial
algebras. Relative to the classical case, many new aspects have to
be taken into account. We describe these groups explicitly when
the underlying polytope is $2$-dimensional. Already this
low-dimensional case provides interesting classes of groups.
\end{abstract}

\maketitle

\section{Introduction}\label{INTR}

Polyhedral $K$-groups are associated to groups of graded
automorphisms of algebras $R[P]$ associated with lattice polytopes
$P$ and their arrangements. These are essentially the homogeneous
coordinate rings of projective toric varieties if the ring $R$ of
coefficients is an algebraically closed field. For the special
case of a unit simplex $P$ the algebra $R[P]$ is a polynomial ring
and the group is just the general linear group. Thus polyhedral
$K$-theory contains ordinary (algebraic) $K$-theory as a special
case. However, for more general polytopes many new aspects have to
be taken into account and new $K$-groups appear.

There are several sources of motivation for the
``polyhedrization'' of $K$-theory. They can be described as
follows.

In the series of papers \cite{BrG1}--\cite{BrG4} we have
investigated polytopal algebras $k[P]$ ($k$ a field) as
generalizations of $k$-vector spaces. Recall that for a lattice
polytope $P\subset\RR^n$ the associated polytopal algebra is the
semigroup algebra $k[P]=k[S_P]$ of the (additive) sub-semigroup
$S_P\subset\ZZ^{n+1}$ generated by $\{(x,1)\ |\ x\in\L_P\}$, where
$\L_P=P\cap\ZZ^n$. Equivalently, $k[P]$ is generated over $k$ by
the lattice points of $P$ which are subject to the binomial
relations reflecting the affine dependencies inside $P$. Polytopal
algebras and their graded $k$-homomorphisms ($\deg(x)=1$ for
$x\in\L_P$) define the {\it polytopal linear category} $\Pol(k)$.
The category $\Vect(k)$ of finitely generated $k$-vector spaces is
a full subcategory of $\Pol(k)$: just consider the subclass of
polytopal algebras associated with unimodular simplices (including
the empty one).

In addition to the natural interest in the polytopal
generalization of linear algebra the motivation also comes from
applications to toric geometry -- in \cite{BrG1} we have derived
strengthened versions (for projective toric varieties) of the well
known results on automorphism groups by Demazure \cite{D} and Cox
\cite{Cox}. In \cite{BrG2} this description has been generalized to
arrangements of projective toric varieties.

These results emphasize the analogy between the groups $\GL_n(k)$
and the groups $\ggr.aut(k[P])$ of graded automorphisms of
polytopal algebras $k[P]$. A crucial r\^ole in the description of
$\ggr.aut(k[P])$ is played by {\it elementary automorphisms} of
polytopal algebras which are just ordinary elementary matrices in
the case the underlying polytope is a unimodular simplex. They
made their first appearance (in the setting of toric varieties) in
\cite{D}, the work that initiated the theory of toric varieties in
the 1970s.

The objects studied in our previous papers include the {\it
polytopal Picard group} of $\Pol(k)$, automorphism groups of
polytopal algebras and retractions of polytopal algebras. Retracts
of free modules are projective modules. Therefore the study of
algebra retracts can be considered as a non-linear variant of
studying the group $K_0$ of a ring. The group $K_1$ compares
automorphisms of free modules to the elementary ones, as does our
theorem on the automorphisms in $\Pol(k)$ (Theorem \ref{plg}
below). Therefore the latter is a non-linear analogue of $K_1$,
and it is natural to push the analogy between $\Vect(k)$ and
$\Pol(k)$ further into higher $K$-theory.

Algebraic $K$-theory (for rings) can to some extent be thought of
as a theory of higher syzygies between elementary matrices.
Accordingly, our goal is to develop the theory of higher syzygies
between elementary automorphisms of polytopal algebras.

However, there is yet another source of interest in polyhedral
$K$-theory. The group of automorphisms of polynomial algebras
(over a field) has been a big challenge for researchers for
several decades. The tame generation conjecture asserts that every
automorphism is a composite of linear and triangular
automorphisms. Only the case of two variables has been settled
(Jung, van der Kulk, in the 40s and 50s \cite{J}, \cite{K}). There
was an attempt in the 70s (Bass, Connell and Wright \cite{BaW},
\cite{Con1}, \cite{Con2}, \cite{ConW}) to develop a {\it
non-linear} $K$-theory based on these groups as non-linear
analogues of $\GL_n$, and on retracts of polynomial rings as
analogues of projective modules. But progress in this area seems
to be blocked.

The relationship between polyhedral and non-linear polynomial
$K$-theories looks as follows:
\begin{gather*}
\begin{CD}
\text{arbitrary automorphisms}@>\text{restricted to}>>
\text{graded automorphisms,}
\end{CD}
\\
\begin{CD}
\text{polynomial rings}@>\text{generalized to}>>\text{polytopal
rings} @>>>\text{polyhedral algebras.}
\end{CD}
\end{gather*}
Moreover, one easily observes that the tame generation conjecture
is equivalent to the generalization of Theorem \ref{plg} to not
necessarily bounded lattice polytopes (actually the case of
positive orthants $\Pi=\RR^n_+$, $n\in\NN$ would already suffice,
see Proposition \ref{tame}). The same relationship exists with the
classical conjecture on retracts of polynomial algebras.

The polytopal linear groups interpolate between the linear and
the `infinite  dimensional' groups (see Shafarevich \cite{Sh}).
Thus they provide a possibility for studying similarities of the
`infinite dimensional' objects and the linear ones within a finite
dimensional framework.

The automorphism groups $\ggr.aut(k[P])$ and their
subgroups of elementary automorphisms $\EE_k(P)$ are linear
groups in a natural way. But unlike $\GL_n(k)$ they are
essentially never reductive groups, as documented by the
list in Theorem \ref{sta2gp}, which describes the
corresponding \emph{stable} groups $\EE(R,P)$ when $\dim
P=2$ (and $k=R$). The upper triangular subgroups of the
$\EE(R,P)$ are contributions of the unipotent radicals of
the corresponding unstable groups $\EE_k(P_i)$, $i\in\NN$.
(For $\dim P>2$ a complete classification as in Theorem
\ref{sta2gp} seems to be very difficult.) Further, the
groups $\EE_k(P)$ are very often non-perfect. Even more is
true: it follows from Demazure's work \cite{D} that the
groups $\ggr.aut(k[P])$ in the important special case of a
smooth variety $\Proj(K[P])$ are always semidirect products
of unipotent groups and reductive groups with root system
of type $A_l$. Therefore, the theory developed here is to
some extent complementary to the theory of Chevalley groups
and their universal central extensions. Milnor's definition
of $K_2$ was motivated by Steinberg's fundamental work
\cite{Stb} who considered such groups over fields. This was
further extended to groups over arbitrary commutative rings
by Stein \cite{St} (see also \cite{KaSt}).

One has to reveal the hidden polytopal geometry behind the
Steinberg relations when elementary matrices are generalized to
elementary automorphisms. This leads us to introduce the class of
\emph{balanced} polytopes (Section \ref{BALANCED}). They allow
various $K$-theoretic constructions (further confirmation of which
is provided in \cite{BrG5}).

One of the main features of the polyhedral theory (and a
difference to Chevalley groups) is that one has to work with
\emph{stable} groups, rather than with the groups
$\ggr.aut_R(R[P])$. This is explained by the absence of inductive
homomorphisms between them (Remark \ref{nfiltun}). The recipe for
the definition of
stable groups is encoded in a purely polytopal construction called
\emph{doubling along a facet} (Section \ref{DOUBL}), which is
minimal with respect to several natural properties (Remark
\ref{informal}(a)). The elementary automorphisms extend to the
``doubled'' polytopal algebras and the final outcome is a perfect
group. In the special case of unit simplices one recovers the
familiar sequence
$$
E_2(R)\subset\dots\subset E_n(R)\subset E_{n+1}(R)\subset\cdots,
\qquad*\mapsto
\begin{pmatrix}
*&0\\
0&1
\end{pmatrix}.
$$
Using such {\it doubling spectra} for a commutative ring $R$ and a
lattice polytope $P$ we define the stable group of elementary
automorphisms $\EE(R,P)$ (Section \ref{STABLE}) and the stable
Steinberg group $\St(R,P)$ (Section \ref{SCHUR}). However, the
situation for general polytopes is more complicated than for unit
simplices -- as remarked, $\EE(R,P)$ is not the inductive
limit of the corresponding unstable groups and one has to encode
the inductive process into the underlying polytopes themselves,
rather than into the automorphism groups. Actually, the groups
obtained depend bivariantly on $R$ and $P$ (more precisely, on the
normal fan of $P$; Section \ref{BIFUN}).

The main results of this paper are:
\begin{itemize}
\item[(i)] The existence of the universal central extension of
polytopal groups
$$
1\to K_2(R,P)\to\St(R,P)\to\EE(R,P)\to1
$$
provided $P$ is balanced (Theorem \ref{ucext}). Unlike in the
situation of unit simplices, the fact that $K_2(R,P)$ coincides
with the center of $\St(R,P)$ is no longer trivial. It follows
from the triviality of the center of $\EE(R,P)$ (Theorem
\ref{trivcen}).

\item[(ii)] The visualization of these groups in the first
non-trivial case when $P$ is a balanced polygon (i.~e.\ $\dim
P=2$). We will see that there exist essentially $6$ different
types (Section \ref{POLYG}).
\end{itemize}

Most of the results below generalize to the bigger class of {\it
polyhedral algebras} \cite{BrG3}, which are composed of polytopal
algebras in the same way as Stanley-Reisner rings are composed of
polynomial algebras. However, we do not pursue this level of
generality; we only use the attribute ``polyhedral'' to indicate
the possibility of generalization. The results of this paper are
essential for our study \cite{BrG5} of higher polyhedral
$K$-groups.

\section{Column vectors and elementary automorphisms}\label{ELAT}

In this section we introduce the notion of elementary
automorphisms which is crucial for the sequel. To this end we
recall the relevant part of \cite{BrG1}. The only difference to
\cite{BrG1} is that we will consider arbitrary commutative rings
instead of fields of coefficients. We do not include any arguments
because they are completely parallel to those for fields.
Throughout the paper $R$ denotes an arbitrary commutative ring.

A polytope in $\RR^n$ is always assumed to be finite and convex.
Further, we will only consider lattice polytopes, i.~e.\ those
with vertices in $\ZZ^n$. The objects to be defined below depend
only on the affine structure of the set $\L_P=P\cap\ZZ^n$ of
lattice points of $P$. Therefore we can always assume that $\RR^n$
is the smallest affine space containing $P$ and that $\ZZ^n$ is
the smallest affine sublattice of $\RR^n$ containing $\L_P$. If
necessary, we simply replace $\RR^n$ by the affine hull $\Aff(P)$
of $P$ and $\ZZ^n$ by $z_0+\sum_{z\in\L_P}\ZZ(z-z_0)$ for some
$z_0\in \L_P$.

Under this assumption let $F$ be a facet of $P$ and choose a point
$z_F\in F$. Then the subgroup
$$
F_\ZZ:=(-z_F+\Aff(F))\cap\ZZ^n\subset\ZZ^n
$$
is isomorphic to $\ZZ^{n-1}$. Moreover, there is a unique group
homomorphism $\la F,-\ra:\ZZ^n\to\ZZ$, written as $x\mapsto\la
F,x\ra$, such that $\Ker(\la F,-\ra)=F_\ZZ$, $\Coker(\la
F,-\ra)=0$, and $\la F,-\ra$ attains its minimum $b_F$ on the set
$\L_P$ at the lattice points of $F$.

The homomorphism $\la F,-\ra$ can be naturally extended to
$\RR^n$. On $\Aff(F)$ it has the constant value $b_F\in\ZZ$. With
this notation we obtain the representation
$$
P=\{x\in\RR^n:\la F,x\ra\ge b_F\text{ for all facets }F\},
$$
of $P$ as an intersection of closed halfspaces.

We can use $\la F,-\ra$ to measure the \emph{height} of a lattice
point $x$ above the facet $F$ by setting
$$
\het_F(x)=\la F,x\ra-b_F.
$$
The function $\het_F$ counts the number of hyperplanes between $F$
and $x$ (in the direction of $P$) that are parallel to, but
different from $\Aff(F)$ and pass through lattice points.

The polytopal algebra $R[P]$ is by definition the $R$-algebra
generated by the lattice points of $P$ which are subject to the
binomial relations reflecting the affine dependencies inside $P$.
Equivalently, $R[P]$ is the semigroup ring $R[S_P$] of the
additive subsemigroup $S_P\subset\ZZ^{n+1}$ generated by
$$
\{(x,1)\ |\ x\in\L_P\}\subset\ZZ^{n+1}.
$$
For a polytope $P\subset\RR^n$ the group
$\Gamma_R(P)=\ggr.aut_R(R[P])$
coincides with $\GL_n(R)$ in the case of the unit $(n-1)$-simplex
$$
P=\Delta_{n-1}=\conv\bigl((0,0,\ldots,0),(1,0,\ldots,0),
\ldots,(0,0,\ldots,1)\bigr)\subset\RR^{n-1}.
$$
The groups $\Gamma_R(P)$ are called {\em polytopal linear groups}
-- they consist of the $R$-points of the corresponding
group schemes (which are defined over $\ZZ$, see
\cite[\S2]{BrG4}).

An element $v\in \ZZ^n$, $v\neq 0$, is a {\em column vector} (for
$P$) if there is a facet $F\subset P$ such that $x+v\in P$ for
every lattice point $x\in P\setminus F$. The facet $F$ is called
{\em the base facet} of $v$. While the points $x\in P\cap\ZZ^n$
are identified with $(x,1)\in\ZZ^{n+1}$, the vector $v$ is to be
identified with $(v,0)\in\ZZ^{n+1}$.

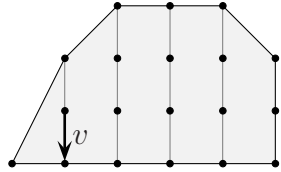
\begin{figure}[htb]
\begin{center}
\begin{pspicture}(0,0)(5,3)
\pspolygon[style=fyp](0,0)(5,0)(5,2)(4,3)(2,3)(1,2)
\psline[linecolor=gray](2,0)(2,3)
\psline[linecolor=gray](3,0)(3,3)
\psline[linecolor=gray](4,0)(4,3)
\psline[linecolor=gray](1,0)(1,2) \multirput(0,0)(1,0){6}{\vertex}
\multirput(1,1)(1,0){5}{\vertex} \multirput(1,2)(1,0){5}{\vertex}
\multirput(2,3)(1,0){3}{\vertex} \rput(1.3,0.5){$v$}
\psline[style=fatline]{->}(1,1)(1,0.05)
\end{pspicture}
\end{center}
\caption{A column structure}\label{FigCol}
\end{figure}

The set of column vectors of $P$ is denoted by $\Col(P)$. A pair
$(P,v)$, $v\in\Col(P)$, is called a {\em column structure} (see
Figure \ref{FigCol}).

The function $\het_F$, defined above on
$\ZZ^n\approx\ZZ^n\times\{1\}$ for each facet $F$ of $P$, has a
unique extension to a linear form on $\RR^{n+1}$. The extension is
again denoted by $\het_F$. For $\het_{P_v}$ we simply write
$\het_v$. It is an easy observation that $x+\het_v(x)\cdot v\in
S_{P_v}\subset S_P$ for any $x\in S_P$.

Column vectors are the dual objects to the roots of the {\it
normal fan} $\mathcal N(P)$ in the sense of Demazure (see
\cite{D}, \cite{Cox}; for normal fans see Section \ref{BIFUN}).

Let $(P,v)$ be a column structure and $\lambda\in R$. Then the
assignment
$$
x\mapsto (1+\lambda v)^{\het_v(x)}x.
$$
gives rise to a graded $R$-algebra automorphism $e_v^\lambda$ of $R[P]$.
(At this point we need the commutativity of $R$.)
Its inverse is $e_v^{-\lambda}$.
Observe that $e_v^\lambda$ becomes an elementary matrix in the special
case when $P=\Delta_{n-1}$, after the identifications
$R[{\Delta_{n-1}}]=R[X_1,\dots,X_n]$ and $\Gamma_R(P)=\GL_n(R)$. Accordingly
$e_v^\lambda$ is called an {\em elementary automorphism}.

The following alternative description of elementary automorphisms
shows that they are always restrictions of elementary matrices if
we interpret the latter as automorphisms of polynomial algebras.

We may assume $\dim P=n$ and $\gp(S_P)=\ZZ^{n+1}$. By a suitable
integral unimodular change of coordinates we may further assume
that $v=(0,-1,0,\dots,0)$ and that $P_v$ lies in the subspace
$\RR^{n-1}$ (thus $P$ is in the upper halfspace). Consider the
standard unimodular simplex $\Delta_n$ (i.~e. the one with
vertices at the origin and the standard coordinate unit vectors).
Clearly, $P$ is contained in a parallel integral shift of
$c\Delta_n$ for a sufficiently large natural number $c$. Then we
have a graded $R$-embedding $R[P]\to R[c\Delta_n]$, the latter
ring being just the $c$th Veronese subring of the polynomial ring
$R[\Delta_n]=R[X_0,X_1,\dots,X_n]$. Moreover, $v=X_0/X_1$. Now the
automorphism of $R[X_0,X_1,\dots,X_n]$ mapping $X_1$ to
$X_1+\lambda X_0$ and leaving all the other variables fixed
induces an automorphism of $R[c\Delta_n]$ which restricts to an
automorphism of $R[P]$. It is nothing but the elementary
automorphism $e_v^{\lambda}$ above.

\begin{proposition}\label{afemb}
Let $v_1,\dots,v_s$ be pairwise different column vectors for $P$ with
the same base facet $F=P_{v_i}$, $i=1,\dots,s$. The mapping
$$
\phi:R^s\to\Gamma_R(P),\qquad
(\lambda_1,\dots,\lambda_s)\mapsto
e_{v_1}^{\lambda_1}\circ\cdots\circ e_{v_s}^{\lambda_s},
$$
is an injective group homomorphism from the additive group $R^s$. In
particular, $e_{v_i}^{\lambda_i}$ and $e_{v_j}^{\lambda_j}$ commute
for all $i,j\in\{1,\dots,s\}$ and $(e_{v_1}^{\lambda_1}\circ\cdots\circ
e_{v_s}^{\lambda_s})^{-1} =e_{v_1}^{-\lambda_1}\circ\cdots\circ
e_{v_s}^{-\lambda_s}.$
\end{proposition}

In what follows $\EE_R(P)$ denotes the subgroup of $\Gamma_R(P)$
generated by the elementary automorphisms. In particular ,
$\EE_R(\Delta_{n-1})=\mathrm{E}_n(R)$ -- the $n$th unstable group
of elementary matrices.

The initial motivation for developing the theory of higher
syzygies between the elementary automorphisms is the main result
of \cite{BrG1} (Theorem \ref{plg} below), which establishes a
complete polytopal generalization of the standard linear algebra
fact that any invertible matrix over a field can be reduced to a
diagonal matrix using elementary transformations on columns (or
rows). Moreover, we have normal forms for such reductions since
the elementary transformations can be carried out in an increasing
order of the column indices.

Let $k$ be a field. Put $d=\dim P+1$. The $d$-torus
$\TT_d=(k^*)^d$ acts naturally on $k[\gp(S_P)]$ via the
substitution
$$
(\xi_1,\dots,\xi_d)(e_i)=\xi_ie_i,\qquad \xi_i\in k^*,\ \
i\in[1,d].
$$
Here $e_i$ is the $i$-th element of a fixed basis of
$\gp(S_P)=\ZZ^d$. Since the action restricts to $k[P]$, one has an
algebraic embedding $\TT_d\subset\Gamma_k(P)$, whose image we
denote by $\TT_k(P)$. It consists precisely of those automorphisms
of $k[P]$ which multiply each monomial by a scalar from $k^*$.

The (finite) automorphism group $\Sigma(P)$ of the semigroup $S_P$ is
also a subgroup of $\Gamma_k(P)$. It is exactly the group of automorphisms
of $P$ as a lattice polytope.

The embedding $\phi$, given by Lemma \ref{afemb}, is an embedding
of algebraic groups over $k$. Denote by $\AA(F)$ the image of
$\phi$. It is an affine space over $k$. Of course, $\AA(F)$ may
consist only of the identity map of $k[P]$ -- namely, if there is
no column vector with base facet $F$.

\begin{theorem}\label{plg}
Let $P$ be a $n$-polytope and $k$ a field.
Every element $\gamma\in\Gamma_k(P)$ has a (not uniquely determined)
presentation
$$
\gamma=\alpha_1\circ\alpha_2\circ\cdots\circ\alpha_r\circ\tau\circ\sigma,
$$
where $\sigma\in\Sigma(P)$, $\tau\in\TT_k(P)$, and
$\alpha_i\in\AA(F_{i})$ such that the facets $F_i$ are pairwise
different and $\#\L_{F_i}\le \#\L_{F_{i+1}}$, $i\in[1,r-1]$.

We have $\dim\Gamma_k(P)=\#\Col(P)+d$ (the left hand side is the Krull
dimension of the affine group scheme $\Gamma_k(P)$), and $\TT_k(P)$ is a
maximal torus in $\Gamma_k(P)$, provided $k$ is infinite.
\end{theorem}

The convex hull in $\RR^n$ of a subset $X\subset\ZZ^n$ will be
called a {\it lattice polyhedron} if it is an intersection of a
{\it finite} family of rational affine halfspaces in $\RR^n$.
(`Rational affine' means that the hyperplane bounding the
halfspace is given by a rational linear equation.) For polyhedra
(in contrast to polytopes) we do no longer require compactness.

Let $P\subset\RR^n$ be a convex lattice polyhedron. As in the case
of polytopes, there exists a height function $\het_F$ for each
facet $F$ of $P$. It is the unique surjective homomorphism
$\ZZ^n\to\ZZ$ that maps $\L_F$ to $0$ and $\L_P$ to $\ZZ_+$. Thus
one can define the notion of column vector for $P$: an element
$v\in\ZZ^n\setminus\{0\}$ is a column vector if there exists a
facet $F\subset P$ such that (i) $x+v\in P$ whenever
$x\in\L_P\setminus F$ and (ii) $\het_F(v)=-1$. Observe that
condition (i) implies (ii) if $P$ is a (bounded) polytope, but
this is not true for (unbounded) polyhedra in general.

One can define the `non-compact' analogues of polytopal algebras
and their elementary automorphisms, using this notion of column
vectors. Clearly, Proposition \ref{afemb} generalizes completely
to the general polyhedral setting.
\bigskip

\noindent {\em The positive orthant and the tameness
conjecture.}\enspace Consider the positive orthant $\Pi=\RR^n_+$.
One observes that $\Gamma_k(\Pi)$ contains a copy of the torus
$(k^*)^{n+1}$, just as in the case of finite polytopes.

We want to relate $k[\Pi]$ to the polynomial ring
$k[X_1,\dots,X_n]$. To this end we set $\ee_0=(0,1)\in S_\Pi$ and
$\ee_i=(e_i,1)\in S_\Pi$, $i\in[1,n]$. Then the assignment
$X_i\mapsto \ee_i$ embeds $k[X_1,\dots,X_n]$ into $k[\Pi]$. (This
is just the embedding coming from the identification of
$\Delta_{n-1}$ with the simplex spanned by $e_1,\dots,e_n$.) We
identify the polynomial ring with its image. Then clearly
$k[\Pi][\ee_0^{-1}]\approx k[X_1,\dots,X_n][[\ee_0^{\pm1}]$, and
every automorphism of the polynomial ring extends to one of $k[\Pi]$
that leaves $\ee_0$ invariant.

Conversely, one has a natural isomorphism $k[\Pi]/(\ee_0-1)\approx
k[X_1,\dots,X_n]$ given by dehomogenization, making
$k[X_1,\dots,X_n]$ a retract of $k[\Pi]$. Therefore every
automorphism of $k[\Pi]$ that leaves $\ee_0$ fixed induces an
automorphism of $k[X_1,\dots,X_n]$.

The analogue of Theorem \ref{plg} for the infinite polyhedron
$\Pi$ is equivalent to the tameness conjecture for
the group of $k$-automorphisms of $k[X_1,\dots,X_n]$, as the next
proposition shows.

Recall that an automorphism $\alpha:k[X_1,\ldots,X_n]\to
k[X_1,\ldots,X_n]$ is called tame if it is a composite of affine
(i.~e.\ linear + constant) automorphisms and triangular
automorphisms, i.~e.\ those of type
$$
X_1\mapsto X_1,\qquad X_i\mapsto X_i+\phi_i,\quad\phi_i\in
k[X_1,\ldots,X_{i-1}],\quad i\in[2,n].
$$
Furthermore let $\TT_k(\Pi)'$ denote the subgroup of toric
automorphisms that leaves $\ee_0$ invariant.

\begin{proposition}\label{tame}
For a field $k$ the subgroup of $\Gamma_k(\Pi)$ generated by
$\EE_k(\Pi)$  and $\TT_k(\Pi)'\approx (k^*)^{n}$ is naturally
isomorphic to the group of tame automorphisms of the polynomial
ring $k[X_1,\ldots,X_n]$. Moreover, under this isomorphism
$\EE_k(\Pi)$ corresponds to the group of tame automorphisms of
Jacobian $1$.
\end{proposition}

The proposition is easily proved (we leave the details to the
reader).

Notice that Proposition \ref{tame} treats the augmented and
non-augmented automorphisms in a uniform way. Moreover, there is
no need to consider the affine and triangular automorphisms
separately (as it is usually done in the literature).

\section{The partial product operation on column vectors}\label{PRODUCTS}

For a polytope $P$ we denote by $\FF(P)$ the set of all facets of
$P$. The following lemma is easily derived from the definition of
column vector and the description of $P$ as the intersection of
closed halfspaces.

\begin{lemma}\label{colv}
A vector $v\in\ZZ^n$ belongs to $\Col(P)$ if and only if there
exists $F\in\FF(P)$ such that
$$
\la F,v\ra=-1\quad\text{and}\quad\la F,v\ra\geq0\ \text{for all
}G\in\FF(P),\ G\neq F.
$$
\end{lemma}

Clearly, this facet $F$ is just the base facet $P_v$.

Now we introduce the notion of a partial product operation on
column vectors and investigate its basic properties. Fix a
polytope $P$ and let $(P,u)$ and $(P,v)$ be column structures on
it.

\begin{definition}\label{prod}
We say that the product $uv$ exists if $u+v\neq 0$ and $x+u\notin
P_v$ for every point $x\in\L_P\setminus P_u$. If $uv$ exists, we
put $uv=u+v$.
\end{definition}

\begin{figure}[htb]
\begin{center}
 \psset{unit=1cm}
 \def\vertex{\pscircle[fillstyle=solid,fillcolor=black]{0.05}}
\begin{pspicture}(-0.3,-0.5)(4,2.5)
 \pspolygon[style=fyp,linecolor=darkgray](0,0)(3,0)(3,2)(2,2)
 \footnotesize
 \multirput(0,1)(1,0){4}{\vertex}
 \multirput(0,0)(1,0){4}{\vertex}
 \multirput(0,2)(1,0){4}{\vertex}
 \psline[style=fatline]{->}(1,1)(0,0)
 \psline[style=fatline]{->}(1,0)(0,0)
 \psline[style=fatline]{->}(2,1)(3,1)
 \psline[style=fatline]{->}(1,1)(1,0)
 \psline[style=fatline]{->}(3,2)(2,1)
 \psline[style=fatline]{->}(3,2)(3,1)
 \rput(-0.2,0.6){$w=uv$}
 \rput(1.3,0.5){$u$}
 \rput(0.6,-0.2){$v$}
 \rput(2.5,0.8){$-v$}
 \rput(2.1,1.5){$w$}
 \rput(4.0,1.5){$u=w(-v)$}
\end{pspicture}
\end{center}
\caption{The product of two column vectors} \label{ProdCol}
\end{figure}
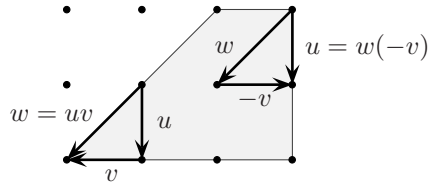

Figure \ref{ProdCol} shows a polytope with all its column
vectors and the two existing products $w=uv$ and $u=w(-v)$.

\begin{proposition}\label{maincrit}
\leavevmode\
\begin{itemize}
\item[(a)]
The following conditions are equivalent:
\begin{enumerate}
\item $uv$ exists;
\item $u+v\in\Col(P)$ and $P_{u+v}=P_u$;
\item $u+v\neq0$ and $\la P_v,u\ra>0$.
\end{enumerate}
In particular, $P_{uv}=P_u$ and $v$ is parallel to $P_u$ (i.~e.\
$\la P_u,v\ra=0$). \item[(b)] $u+v\in\Col(P)$ if and only if
exactly one of the products $uv$ and $vu$ exists. \item[(c)]
Suppose $u,v\in\Col(P)$, $u\neq-v$, and $\la P_v,u\ra>0$. Then
$u+iv\in\Col(P)$ for all $i\in[1,\la P_v,u\ra]$ and
$u+iv\notin\Col(P)$ for all $i>\la P_v,u\ra$. \item[(d)] The
following are equivalent for $v\in\Col(P)$:
\begin{enumerate}
\item $-v\in\Col(P)$; \item there exist facets $F,G$ of $P$ such
that $\la F,v\ra=-1$, $\la G,v\ra=1$, and $\la H,v\ra=0$ for every
facet $H\neq F,G$; \item there exists $w\in\Col(P)$ and facets
$F,G$ such that $\la F,v\ra=-1$, $\la G,v\ra=1$, $\la F,w\ra=1$,
$\la G,w\ra=-1$.
\end{enumerate}
\item[(e)] If $w=uv$ and $-w\in\Col(P)$, then $-u,-v\in\Col(P)$ as
well.
\end{itemize}
\end{proposition}

\begin{proof} (a) The equivalence of (1) and (2) as well as
the implication (3)$\implies$(2) are straightforward. Now assume
$uv$ exists and $\la P_v,u\ra\leq0$. If $\la P_v,u\ra<0$, then
$P_u=P_v$ according to Lemma \ref{colv} and therefore $\la
F_u,uv\ra=-2$, which is impossible (again by Lemma \ref{colv}). If
$\la P_v,u\ra=0$ then $P_u\neq  P_v$ and, hence, there is an
element $x\in\L_{P_v}\setminus P_u$. But then the point $x+u+v$ is
outside $P$ -- a contradiction.

(b) Assume $u+v\in\Col(P)$. According to (a) we have to show that
exactly one of the inequalities $\la P_v,u\ra>0$ and $\la
P_u,v\ra>0$ holds. The opposite (strict) inequalities are excluded
by the same reasons as in the proof of (a).

Next we exclude that simultaneously $\la P_v,u\ra=0$ and $\la
P_u,v\ra=0$. Assume that $\la P_v,u\ra=0$ and $\la P_u,v\ra=0$. In
particular $P_u\neq P_v$ and, therefore, one of these base facets
differs from $P_{u+v}$. If $P_u\neq P_{u+v}$, then there is a
point $x\in\L_{P_u}\setminus P_{u+v}$ and by the equality $\la
P_u,v\ra=0$ we get a contradiction because $\het_u(x+u+v)=-1$,
that is, $x+u+v$ is outside $P$.

It only remains to show that  $\la P_v,u\ra>0$ and $\la
P_u,v\ra>0$ can not hold simultaneously. But if this were the
case, then for every lattice point $x\in\L_P\setminus P_u$ we
would have
$$
x,\ x+u,\ x+u+v,\ x+u+v+u,\ x+u+v+u+v,\ldots\in\L_P.
$$
This is a contradiction, because $x+c(u+v)\notin P$ for $c\gg0$.

(c) The proof is straightforward.

(d) Suppose $-v\in\Col(P)$. Then condition (2) is satisfied for
$F=P_v$, $G=P_{-v}$: note that $\la H,v-v\ra=0$ and $\la
H,v\ra,\la H,-v\ra\ge0$ for $H\neq F,G$. Similarly (3) holds with
$w=-v$.

If condition (2) holds, then $-v$ is clearly a column vector with
base facet $G$: $x+(-v)\in P$ for all $x\in\L_P\setminus G$.

If (3) holds, then $\la H,v+w\ra\ge 0$ for all facets $H$. The
only vector satisfying this condition is $0$.

(e) Let $F=P_w=P_u$, $G=P_{-w}$. One has $\la F,v\ra=0$, and there
exists a facet $E$ with $\la E,v\ra>0$. If $E\neq G$, then $\la
E,w\ra =\la E, u+v\ra>0$, and this is impossible for $E\neq G$.
Thus $E=G$, and since $\la G,w\ra=1$, it follows that $\la
G,u\ra=0$, $\la G,v\ra=1$.

Now consider the facet $P_v$. As seen already, it is different
from $F$ and $G$. So $\la P_v,w\ra=0$ forces $\la P_v,u\ra=1$.
Moreover, for all facets $H\neq F,G,P_v$ we must have $\la
H,u\ra=\la H,v\ra=0$. In view of (d.2) both $u$ and $v$ are
invertible column vectors.
\end{proof}

\begin{corollary}\label{assoc}
Let $u,v,w\in\Col(P)$.
\begin{itemize}
\item[(a)]
If both $uv$ and $vw$ exist and $u+v+w\neq 0$ then the products
$(uv)w$ and $u(vw)$ also exist and, clearly, $(uv)w=u(vw)$.
\item[(b)]
If $vw$ and $u(vw)$ exist and $u+v\not=0$ then $uv$ also exists,
while the existence of $uv$ and $(uv)w$ does in general not imply
the existence of $vw$, even if $v+w\not=0$
\end{itemize}
\end{corollary}

\begin{figure}[htb]
\begin{center}
\psset{unit=2cm}
\def\vertex{\pscircle[fillstyle=solid,fillcolor=black]{0.05}}
\begin{pspicture}(-0.5,-0.7)(1,0.7)
\psset{viewpoint=-1 -2 -1.5} \ThreeDput[normal=0 0 1](0,0,0){
  \pspolygon[linewidth=0pt, style=fyp, fillcolor=medium](1,0)(0,0)(0,1)(1,1)
} \ThreeDput[normal=0 -1 0]{
  \pspolygon[linewidth=0pt, style=fyp](1,0)(0,0)(0,1)
}
\ThreeDput[normal=1 0 0]{
  \pspolygon[linewidth=0pt, style=fyp](1,0)(0,0)(0,1)
}
\ThreeDput[normal=0 0 1](0,0,0){
  \rput(0,0){\vertex}
  \rput(0,1){\vertex}
  \rput(1,0){\vertex}
  \rput(1,1){\vertex}
  \psline[linestyle=dashed](0,1)(0,0)
  \psline[linewidth=1.5pt]{->}(0,0)(1,0)
  \psline[linewidth=1.5pt]{->}(1,0)(1,1)
  \psline(0,1)(1,1)
}
\ThreeDput[normal=0 -1 0]{
  \psline(0,1)(1,0)
  \psline[linewidth=1.5pt]{->}(0,1)(0,0)
}
 \ThreeDput[normal=1 0 0]{\psline(0,1)(1,0)}
 \ThreeDput[normal=1 -1 0]{\psline(0,1)(1.41,0)}
 \ThreeDput(0,0,1){\rput(0,0){\vertex}}
 \rput(-0.12,0.1){\footnotesize$u$}
 \rput(0.6,-0.05){\footnotesize$v$}
 \rput(0.75,-0.7){\footnotesize$w$}
\end{pspicture}
\end{center}
\caption{The pyramid over the unit square} \label{ThreeCol}
\end{figure}
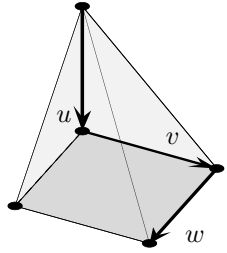

\begin{proof}
(a) That $(uv)w$ exists follows from the definition of the
product and the equality $P_{uv}=P_u$. Thus, by Proposition
\ref{maincrit}(a) we only need to
show that $\la P_{vw},u\ra>0$. But the same proposition implies
$\la P_{vw},u\ra=\la P_v,u\ra>0$.

(b) That the existence of $vw$ and $u(vw)$ implies that of $uv$
follows from $\la P_v,u\ra=\la P_{vw},u\ra>0$.

The example (see Figure \ref{ThreeCol})
\begin{gather*}
P=\conv\{(0,0,0),(0,0,1),(1,0,0),(0,1,0),(1,1,0)\}\subset\RR^3,\\
u=(0,0,-1),\ v=(1,0,0),\ w=(0,1,0)\in\Col(P),
\end{gather*}
completes the proof because $uv$ and $(uv)w$ exist and $vw$ does
not exist.
\end{proof}

\section{Doubling along a facet}\label{DOUBL}

We now introduce an operator on the set of lattice polytopes,
called \emph{doubling along a facet}.

Let $P\subset\RR^n$ be a polytope and $F\subset P$ be a facet. For
simplicity we assume that $0\in F$, a condition that can be
satisfied by a parallel translation of $P$. Denote by
$H\subset\RR^{n+1}$ the $n$-dimensional linear subspace that
contains $F$ and whose normal vector is perpendicular to that of
$\RR^n=\RR^n\oplus0\subset\RR^{n+1}$ (with respect to the standard
scalar product on $\RR^{n+1}$). Then the upper half space
$H\cap\bigl(\RR^n\times\RR_+\bigr)$ contains a congruent copy of
$P$ which differs from $P$ by a $90^\circ$ rotation. Denote the
copy by $P^{\vt_F}$, or just by $P^\vt $ if there is no danger of
confusion.

Note that $P^\vt$ is not always a lattice polytope with respect to
the standard lattice $\ZZ^{n+1}$. However, it is so with respect
to the sublattice $(\ZZ^n)^{\vt_F}$ which is the image of $\ZZ^n$
under the $90^\circ$ rotation.

The operator of doubling along a facet is then defined by
$$
P^{\sq_F}=\conv(P,P^\vt )\subset \RR^{n+1}.
$$
For typographical reasons we will sometimes write $\sqd(P,F)$
instead of $P^{\sq_F}$.

The doubled polytope is a lattice polytope with respect to the
subgroup $(\ZZ^n)^{\sq_F}=\ZZ^n+(\ZZ^n)^{\vt_F}$. After a change
of basis in $\RR^{n+1}$ that does not affect $\RR^n$ we can
replace $(\ZZ^n)^{\sq_F}$ by $\ZZ^{n+1}$, and consider $P^{\sq_F}$
as an ordinary lattice polytope in $\RR^{n+1}$. In what follows,
whenever we double a lattice polytope $P\subset\RR^n$ along a
facet $F$, the lattice of reference in $\RR^{n+1}$ is always
$\ZZ^n+(\ZZ^n)^{\vt_F}$. For simplicity of notation this lattice
will be denoted by $\ZZ^{n+1}$. (We are grateful to the referee
for pointing out an incorrectness in the previous version of this
article.)

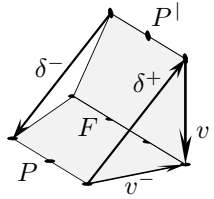
\begin{figure}[htb]
\begin{center}
\psset{unit=1cm}
\def\vertex{\pscircle[fillstyle=solid,fillcolor=black]{0.07}}
\begin{pspicture}(-2.8,0)(1,2)
\psset{viewpoint=3.5 2 -1.5}
\footnotesize
\ThreeDput[normal=0 0 -1](0,0,0){
  \pspolygon[style=fyp](0,0)(3,0)(2,1.5)(0,1.5)(0,0)
  \multirput(0,0)(1,0){4}{\vertex}
  \multirput(0,1.5)(1,0){3}{\vertex}
  \psline[linewidth=1.5pt]{->}(0,1.5)(0,0) 
 }
\ThreeDput[normal=0 -1 0](0,0,0){
  \pspolygon[style=fyp](0,0)(3,0)(2,1.5)(0,1.5)(0,0)
  \multirput(0,1.5)(1,0){3}{\vertex}
  \psline[linewidth=2.0pt]{->}(0,1.5)(0,0) 
 }
\ThreeDput[normal=-1 0 0](0,0,0){
  \psline[linewidth=1.0pt]{->}(1.5,0)(0,1.5)
 }
\ThreeDput[normal=-1 0 0](2,0,0){
  \psline[linewidth=1.0pt]{->}(0,1.5)(1.5,0)
 }
 \rput(-2.1,-0.1){$P$}
 \rput(-0.25,2.0){$P^\vt $}
 \rput(-1.3,0.5){$F$}
 \rput(-1.8,1.3){$\delta^-$}
 \rput(-0.5,1.1){$\delta^+$}
 \rput(-0.6,-0.25){$v^-$}
 \rput(0.3,0.5){$v^\vt $}
\end{pspicture}
\caption{Doubling along the facet $F$} \label{Doubling}
\end{center}
\end{figure}

The construction is illustrated by Figure \ref{Doubling}.

Sometimes $P$ will be referred to by $P^{-_F}$, or just by $P^-$.

In the special case when $F=P_v$ for some $v\in\Col(P)$ we will
use the self explanatory notations $P^{\sq_v}$, $\sqd(P,v)$,
$P^{-_v}$, $P^{|_v}$, $v^-$, $v^\vt $, $\Col(P)^-$ and
$\Col(P)^\vt $.

It is clear from the definition that $\dim P^{\sq_F}=\dim P+1$ and
that the correspondence
$$
\Psi_P:\FF(P)\setminus\{F\}\to\FF(P^{\sq_F}),\qquad\Psi(G)=\conv(G^-,G^\vt
),
$$
extends to a bijection
$$
\bigl(\FF(P)\setminus\{F\}{\big
)}\cup\{P^-,F\}\approx\FF(P^{\sq_F}),
$$
for which $P^-\mapsto P^-$ and $F\mapsto P^\vt $. (Here $G^\vt $
denotes the facet of $P^\vt $ that corresponds to $G=G^-$.) This
bijection will again be denoted by $\Psi_P$.

The polytope $P^\sq$ has two distinguished column vectors
$\delta_P^+$ and $\delta_P^-$, which are the lattice unit vectors
in $\ZZ^{n+1}$ parallel to the lines connecting the points
$x^-\in\L_{P^-}\setminus F$ with the corresponding points $x^\vt
\in\L_{P^\vt }$. We choose $\delta_P^+$ in such a way that
$(P^\sq)_{\delta_P^+}=P^\vt $ and $(P^\sq)_{\delta_P^-}=P^-$.
Clearly $\delta_P^-=-\delta_P^+$. For formal correctness we should
write $\delta_{F,P}^+$ and $\delta_{F,P}^-$ or $\delta^+(F,P)$ and
$\delta^-(F,P)$, (as we will do in Sections \ref{STABLE} and
\ref{SCHUR}). But in most of the cases it will be clear from the
context which facet has been used for the doubling, and we will
simply write $\delta_P^+$ and $\delta_P^-$, or just $\delta^+$ and
$\delta^-$.

In almost all cases doublings will be carried out along base
facets of column vectors $v$. Then we simplify the notation
$\delta^+_{P_v,P}$ to $\delta^+_{v,P}$ etc.

The following equations are easily observed:
\begin{gather*}
\la\Psi(G),\delta^+\ra=\la\Psi(G),\delta^-\ra=0 \quad \text{for
all}\
G\in\FF(P)\setminus\{F\},\tag{$\ref{DOUBL}_1$}\\
\la P^-,\delta^+\ra=\la P^\vt ,\delta^-\ra=1, \tag{$\ref{DOUBL}_2$}\\
\la G,z\ra=\la\Psi(G),z\ra \quad \text{for all}\ z\in\ZZ^n,\
G\in\FF(P).\tag{$\ref{DOUBL}_3$}
\end{gather*}
(In $(\ref{DOUBL}_3)$ the pairings are considered for $P$ and $P^\sq$
respectively and
$\ZZ^n$ is thought of as the subgroup $\ZZ^n\oplus0\subset\ZZ^{n+1}$.)

In view of Lemma \ref{colv} the equations $(\ref{DOUBL}_3)$ imply

\begin{lemma}\label{colext}
Let $F\subset P$ be a facet and $v\in\Col(P)$ be a column vector.
Then $v\in\Col(P^{\sq_F})$.
\end{lemma}

Of crucial importance for what follows is

\begin{corollary}\label{decomp}
For any $v\in\Col(P)$ the following equations hold in
$\Col(P^{\sq_v})$:
$$
v=v^-=\delta^+v^\vt ,\qquad v^\vt =\delta^-v^-.
$$
\end{corollary}

Let $\{F_1,\ldots,F_m\}\subset\FF(P)$ be a system of facets. Then
we get the system of polytopes $(P_0,P_1,\ldots,P_m)$ and facets
$G_{ij}$  defined recursively as follows
\begin{alignat*}2
P_0&=P,\qquad& G_{0j}&=F_j,\ j\in [1,m]\\
P_i&=P_{i-1}^{\sq_{G_{i-1,i}}}, \qquad&
G_{ij}&=\Psi_{G_{i-1i}}(G_{i-1j}),\ j\in [1,m],\ i\in [2,m].
\end{alignat*}

We will use the notation
$$
P_m=P^{\sq_{(F_1,\ldots,F_m)}}.
$$
In other words, $P^{\sq_{(F_1,\ldots,F_m)}}$ is the polytope we
get after $m$ successive doublings, starting from $P$, along the
facets corresponding to $F_1,\dots,F_m$.

For any permutation $\sigma\in\Sigma_m$ we can form the polytope
$$
P^{\sq_{(F_{\sigma(1)},\ldots,F_{\sigma(m)})}}.
$$
The point is that the resulting polytope is independent of
$\sigma$. This will be important in our definition of polyhedral
$K$-groups:

\begin{proposition}\label{welldef}
The polytopes $P^{\sq_{(F_1,\ldots,F_m)}}$ and
$P^{\sq_{(F_{\sigma(1)},\ldots,F_{\sigma(m)})}}$ are naturally
isomorphic as lattice polytopes.
\end{proposition}

\begin{proof}
We  define the mapping
$$
\Theta:P^{\sq_{(F_1,\ldots,F_m)}}\to
P^{\sq_{(F_{\sigma(1)},\ldots,F_{\sigma(m)})}}
$$
as follows. For every lattice point $x\in
P^{\sq_{(F_1,\ldots,F_m)}}$ consider the recursively given
sequence of nonnegative integers:
\begin{multline*}
 h_m=\la P_{m-1},x\ra,\
 h_{m-1}=\la P_{m-2},x+h_m\delta^-_{P_m}\ra,\dots,\\
 h_1=\la P_0,x+h_m\delta^-_{P_m}+\cdots+h_3\delta^-_{P_3}+
    h_2\delta^-_{P_2}\ra.
\end{multline*}
We also have the lattice point
$$
y=x+\sum_{i=1}^m h_i\delta^-_{P_i}\in P_0=P.
$$
To the data $(h_m,\ldots,h_1)$ and $y$ we now associate the lattice
point
$$
\Theta(x)=y+\sum_{i=1}^m h_{\sigma(i)}\delta^+_{P^*_i}\in P^*_m,
$$
where the sequence $(P_0^*,P_1^*,\ldots,P_m^*)$ is related to
$P^{\sq_{(F_{\sigma(1)},\ldots,F_{\sigma(m)})}}$ in the same way
as $(P_0,P_1,\ldots,P_m)$ to $P^{\sq_{(F_1,\ldots,F_m)}}$. It is
straightforward to see that we do not get outside the source polytope
and that the mapping $\Theta$ is an affine isomorphism respecting
the lattice structures.
\end{proof}

\begin{remark}\label{informal}
(a) Informally, the polytope $P^{\sq_v}$, $v\in\Col(P)$ is the
universal solution (i.~e. minimal and applicable to all polytopes)
to the problem posed by the following 3 conditions:
\begin{itemize}
\item[(1)]
$P\subset P^{\sq_v}$ is a facet,
\item[(2)]
$\Col(P)\subset\Col(P^{\sq_v})$,
\item[(3)]
$v$ is decomposed into a product of two column vectors in $P^{\sq_v}$.
\end{itemize}
These three properties are crucial for defining the stable group
of elementary automorphisms $\EE(R,P)$ which will turn out to be
perfect (Section \ref{STABLE}).

(b) In the special case of an algebraically closed field $R=k$ the
ring $k[P^{\sq_F}]$ admits an algebro-geometric characterization
(for simplicity we assume $P$ is normal, i.~e.\ $S_P$ is a normal
semigroup): $\Spec(k[P^{\sq_F}])$ is the normalization of $X$ in the
pull-back diagram of toric varieties
$$
\begin{CD}
X@>>>\Spec(k[P])\\
@VVV @VVV\\
\Spec(k[P])@>>>\Spec(k[F]),
\end{CD}
$$
where $\Spec(k[P])\to\Spec(k[F])$ is induced by the identity embedding
$k[F]\to k[P]$. It is a split embedding whose left inverse is the
$k$-homomorphism
$$
\pi:k[P]\to k[F],\qquad\pi(x)=0\ \text{if}\ x\in\L_P\setminus F\
\text{and}\ \pi(x)=x\ \text{if}\ x\in F.
$$
The resulting closed embedding $\Spec(k[F])\subset\Spec(k[P])$
is an equivariant divisor with respect to the big torus
$(k^*)^{\dim P+1}\subset k[P]$. Conversely, any equivariant divisor
is always associated to a facet of $P$.
\end{remark}

\section{Balanced polytopes}\label{BALANCED}

A polytope $P$ is called {\it balanced} if $\la F_u,v\ra\leq 1$
for all $u,v\in\Col(P)$. By Lemma \ref{colv} $P$ is balanced if
and only if $|\la F_u,v\ra|\leq 1$ for all $u,v\in\Col(P)$.

The next lemma shows that column vectors behave well with doublings along
facets.

\begin{lemma}\label{colbal}
Let $P$ be a balanced polytope and $F$ be one of its facets. Then
$$
\Col(P^{\sq_F})=\Col(P)^-\cup\Col(P)^\vt
\cup\{\delta^+,\delta^-\}.
$$
\end{lemma}

\begin{proof}
Pick $v\in\Col(P^{\sq_F})$. First notice that it is impossible to
have simultaneously $\la P^\vt ,v\ra<0$, $\la P^-,v\ra>0$ and
$v\neq \delta^+$. In fact, if these inequalities held, then the
product $\delta^-v$ existed because $P^\vt $ would be the base
facet for $v$, and Proposition \ref{maincrit}(a) applies. But
similar arguments show that the product $v\delta^-$ existed as
well -- a contradiction with Proposition \ref{maincrit}(b).

By symmetry the inequalities $\la P^\vt ,v\ra>0$, $\la P^-,v\ra<0$
for $v\neq  \delta^-$ are also impossible. It only remains to
exclude the case in which $\la P^\vt ,v\ra>0$ and $\la
P^-,v\ra>0$. Then we can conclude that $v$ is parallel to either
$P^\vt $ or $P^-$, as claimed.

But if $\la P^\vt ,v\ra>0$ and $\la P^-,v\ra>0$ then the product
$v\delta^-$ exists and, simultaneously,
$$
\la P^\vt ,v\delta^-\ra= \la P^\vt ,v\ra+\la P^\vt
,\delta^-\ra\geq2
$$
which is impossible because $P$ is balanced.
\end{proof}

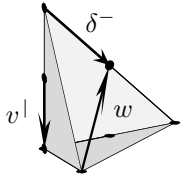
\begin{figure}[htb]
\begin{center}
\psset{unit=1cm}
\def\vertex{\pscircle[fillstyle=solid,fillcolor=black]{0.07}}
\begin{pspicture}(0,-0.3)(1,2)
\psset{viewpoint=-3.5 -2 -1.5} \footnotesize
\ThreeDput[normal=-1 0 0](0,0,0){
  \pspolygon[style=fyp](0,0)(0,2)(2,0)
  }
\ThreeDput[normal=0 0 -1](0,0,0){
  \pspolygon[style=fyp, fillcolor=medium](0,0)(0,2)(1,0)
  \multirput(0,0)(1,0){2}{\vertex}
  \multirput(0,0)(0,1){3}{\vertex}
  \psline[linewidth=1pt]{->}(0,1)(0,0)
 }
\ThreeDput[normal=0 -1 0](0,0,0){
  \pspolygon[style=fyp, fillcolor=medium](0,0)(0,2)(1,0)
  \multirput(0,0)(0,1){3}{\vertex}
  \psline[linewidth=2pt]{->}(0,1)(0,0) 
 }
\ThreeDput[normal=-1 0 0](0,0,0){
  \psline(2,0)(0,2)
  \psline[linewidth=1pt]{->}(0,2)(1,1) 
  \rput(1,1){\vertex}
 }
\ThreeDput[normal=-1 1 0](1,0,0){
  \psline[linewidth=3pt]{->}(0,0)(1.4,1) 
 }
 \rput(1.05,0.5){$w$}
 \rput(0.75,1.7){$\delta^-$}
 \rput(-0.35,0.5){$v^\vt $}
\end{pspicture}
\caption{A new column vector} \label{NewCol}
\end{center}
\end{figure}

\begin{remark}\label{newcol}
(a) Clearly, the union in Lemma \ref{colbal} may not be disjoint
as
$$
\Col(P)^-\cap\Col(P)^\vt =\{v\in\Col(P)\ |\ \la F,v\ra=0\}.
$$

(b) The essential data which determine the partial product
structure of column vectors are the heights $\la F,v\ra$
where $F$ runs through the facets and $v$ through the
column vectors of $P$. Therefore it is useful to clarify
how these data change under doubling. Using the equations
($\ref{DOUBL}_1$), ($\ref{DOUBL}_2$) and ($\ref{DOUBL}_3$),
and the lemma above we can describe the product structure
on $\Col(P^{\sq_F})$: for $u,v\in\Col(P)$ one has
\begin{gather*}
w=uv\iff w^-=u^-v^-\iff w^\vt=u^\vt v^\vt,\\
F=P_v\implies v^-=\delta^+v^\vt\ \text{and}\ v^\vt=\delta^-v^-,\\
-v\in\Col(P)\ \text{and}\ F=P_v\implies v^-(-v)^\vt=\delta^+\ \text{and}\ v^\vt(-v^-)=\delta^-,
\end{gather*}
and there exist no other products in $\Col(P^{\sq_F})$.

(c) If $P$ is not balanced then $P^{\sq_F}$ can have essential new
column vectors, even in the special case when $F=P_v$ for some
$v\in\Col(P)$. For instance, consider the triangle
$\conv((-2,0),(0,0),(0,1))$, which is not balanced, and its column
vector $v=(1,0)$. Then the column vector
$w=(-1,-1,1)\in\Col(P^{\sq_v})$ is not of the type mentioned in
Lemma \ref{colbal}; see Figure \ref{NewCol}.
\end{remark}

\begin{corollary}\label{balpers}
If $P$ is a balanced polytope and $F\subset P$ a facet, then
$P^{\sq_F}$ is balanced as well.
\end{corollary}

This follows immediately from the equations $(\ref{DOUBL}_1)$ and
$(\ref{DOUBL}_2)$ and Lemma \ref{colbal}.

Corollary \ref{balpers} shows that balanced polytopes exist in
abundance. One just starts from a single such polytope and doubles
it successively along arbitrary facets. In \S\ref{POLYG} we will
see that there are infinite families of balanced polytopes even in
the plane.

\section{Commutators of elementary automorphisms}\label{COMMUTATORS}

The main result of this section is the following

\begin{theorem}\label{comrel}
Let $R$ be a ring, $P$ be a polytope, $\lambda,\mu\in R$
and $u,v\in\Col(P)$. Assume $u+v\neq 0$. Then for the commutator
of the automorphisms $e_u^\lambda,e_v^\mu\in\EE_R(P)$ one has
$$
[e_u^\lambda,e_v^\mu]=
\begin{cases}
\prod_{i=1}^n e_{u+iv}^{-\binom{n}{i}\lambda\mu^i}&
\text{if}\ uv\ \text{exists, }n=\la P_v,u\ra,\\
\\
\prod_{i=1}^n e_{v+iu}^{\binom{n}{i}\mu\lambda^i}&
\text{if}\ vu\ \text{exists, }n=\la P_u,v\ra,\\
\\
1&\text{if}\ u+v\notin\Col(P).
\end{cases}
$$
\end{theorem}

The expressions in the first two equations make sense due to
Proposition \ref{maincrit}(c).

\begin{corollary}\label{balcomrel}
Assume $P$ is a balanced polytope and $u,v\in\Col(P)$,
$u+v\neq 0$. Then for all $\lambda,\mu\in R$ we have
$$
[e_u^\lambda,e_v^\mu]=
\begin{cases}e_{uv}^{-\lambda\mu}&\text{if}\ uv\ \text{exists},\\
e_{vu}^{\mu\lambda}&\text{if}\ vu\ \text{exists},\\
1&\text{if}\ u+v\notin\Col(P).
\end{cases}
$$
\end{corollary}

\begin{remark}\label{Steinberg}
Corollary \ref{balcomrel} is the generalization of Steinberg's
relations between elementary matrices to balanced polytopes. In
order to find the classical Steinberg relation $[e_{ij}^\lambda,
e_{jk}^\mu]=e_{ik}^{\lambda\mu}$ in the corollary one must observe
that in our setting the configuration $e_{ij}e_{jk}$ corresponds
to the existence of $vu$ if we associate with $e_{ij}$ the column
vector $v_i-v_j$ where $v_1,\dots,v_n$ are the vectors of the
canonical basis of $R^n$, and simultaneously the vertices of the
unit $n$-simplex.

That we associate $v_i-v_j$ with $e_{ij}$ (and not $v_j-v_i$) is
forced by the our notation in which we add column vectors on the
right. Thus the successive addition of first $u$ and then $v$
corresponds to the product $uv$.
\end{remark}

\begin{proof}[Proof of Theorem \ref{comrel}]
First consider the case when $uv$ exists. As a lattice polytope
$P$ decomposes into lattice polytopal layers parallel to the
affine plane spanned by $u$, $v$ and $uv$, i.~e. we consider the
{\it maximal lattice} polytopes in the sections of $P$ with
2-planes parallel to $\RR u+\RR v$.

Clearly, the automorphisms $\epsilon_1=[e_u^\lambda,e_v^\mu]$ and
$\epsilon_2=e_{u+iv}^{-\binom{n}{i}\lambda\mu^i}$ restrict to
automorphisms of $R[P']$ for each of these layers $P'$ and we have
to check that $\epsilon_1$ and $\epsilon_2$ are the same layer by
layer.

If a layer $P'\subset P$ has dimension $<2$, then at least two of
the vectors $u$, $v$ and $uv$ do not belong to $\Col(P')$. On the
other hand, if $u\in\Col(P')$ or $uv\in\Col(P')$, then one easily
deduces from
Definition \ref{prod} that all the three vectors $u$, $v$
and $uv$ belong to $\Col(P')$. Therefore, $u,uv\notin\Col(P')$.
But this means that $P'\subset P_u=P_{u+iv}$, $i\in[1,n]$.
In particular, the automorphisms $e_u^{\lambda}$ and
$e_{u+iv}^{-\binom{n}{i}\lambda\mu^i}$
restrict to the identity automorphism of $R[P']$ and the
commutator equality becomes trivial on the layer $P'$.

So we only need to consider the 2-dimensional layers $P'$. Then
$u,v,uv\in\Col(P')$ and $uv$ exists as a product in $\Col(P')$.

Since $\gp(S_P)=\ZZ v\oplus\gp(S_{P_v})$ we have
$\la P'_v,u\ra=\la P_v,u\ra$. Without loss of generality we can
therefore assume $P\subset\RR^2$ and $\gp(S_P)=\ZZ^3$.

By a suitable affine integral unimodular change of the coordinates
in $\RR^2$ one can achieve that $P\subset\RR\times\RR_+$,
$P\cap(\RR,0)$ is an edge of $P$, and $u=(0,-1)$, $v=(1,0)$.
Moreover, by a parallel shift we can also assume that $(\la
P_v,u\ra,0)$ is the lowest right vertex of $P$.

Consider the triangle
$$
\Delta=\conv((0,0),(\la P_v,u\ra,0),(0,1))\subset P
$$
(see Figure \ref{EssTri}).
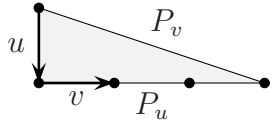
\begin{figure}[htb]
\begin{center}
\psset{unit=1cm}
\def\vertex{\pscircle[fillstyle=solid,fillcolor=black]{0.07}}
\begin{pspicture}(1,0)(2,1)
\pspolygon[style=fyp](0,0)(3,0)(0,1)
\multirput(0,0)(1,0){4}{\vertex} \rput(0,1){\vertex}
\psline[linewidth=1pt]{->}(0,1)(0,0)
\psline[linewidth=1pt]{->}(0,0)(1,0) \rput(1.7,0.8){$P_v$}
\rput(-0.3,0.5){$u$} \rput(1.5,-0.3){$P_u$} \rput(0.5,-0.2){$v$}
\end{pspicture}
\caption{The essential triangle} \label{EssTri}
\end{center}
\end{figure}
We have $\gp(S_P)=\gp(S_\Delta)$. On the other hand
$L_{P_u}\subset k[P]$ is point-wise fixed under automorphisms of
type $e_{u+iv}^*$, $i\in[0,\la P_v,u\ra]$ and $[e_u^*,e_v^*]$
($*\in R$). Therefore we only need to check that the two
automorphisms
$$
[e_u^\lambda,e_v^\mu] \quad\text{and}\quad\prod_{i=1}^n
e_{u+iv}^{-\lambda\binom{n}{i}\mu^i}\in\EE_R(P)
$$
coincide at the vertex $(0,1)$ of $\Delta$.

Observe that all the automorphisms involved restrict to elements
of $\EE_R(\Delta)$. Therefore, after the natural identification
$$
k[\Delta]=k[X^n,X^{n-1}Y,\ldots,XY^{n-1},Y^n][Z],\qquad n=\la
P_v,u\ra,
$$
in such a way that $e_{u+iv}^\lambda(Z)=Z+\lambda X^{n-i}Y^i$ and
$e_v^\lambda(X^{n-i}Y^i)=(X+\lambda Y)^{n-i}Y^i$ for $i\in[0,n]$
the problem has been reduced to the equation
$$
[e_u^\lambda,e_v^\mu](Z)=\prod_{i=1}^n
e_{u+iv}^{-\lambda\binom{n}{i}\mu^i}(Z)
$$
the verification of which is accomplished by a routine
computation.

The case when $vu$ exists is just a consequence of the previous
case.

Now assume $u+v\notin\Col(P)$. In this situation, too, we could
check the desired equality layer by layer. In fact, using
Proposition \ref{maincrit}(a) it is easy to show that
$u+v\notin\Col(P')$ for any of the layers $P'$ as above. But the
arguments below do not get simplified by the consideration of
layers.

By Proposition \ref{maincrit}(c) none of the products $uv$ and
$vu$ exists. By Lemma \ref{colv} and Proposition \ref{maincrit}(a)
there are only two possibilities: either $P_u=P_v$ and the claim
follows from Proposition \ref{afemb} or $\la P_v,u\ra=\la
P_u,v\ra=0$, what we will assume.

Since $v$ is parallel to $P_u$, we get $v\in\Col(P_u)$. Then by
Lemma \ref{colv} there is exactly one facet $g\subset P_u$ such
that $\la g,v\ra=-1$  and for all the other facets $g'\subset P_u$
we have $\la g,v\ra\geq0$ (the pairing is considered on $P_u$).
Now $g$ extends to a unique facet $G\subset P$, $G\neq P_u$, and we see that
$\la G,v\ra<0$. Again by Lemma \ref{colv} $G=P_v$. Using similar
arguments for $u$ and $P_v$ we arrive at the conclusion that the
two facets $P_u,P_v\subset P$ meet along a common facet (a
codimension 2 face of $P$) and, moreover, $P_u\cap P_v$ is a base
facet both for $u\in\Col(P_v)$ and $v\in\Col(P_u)$.

Assume for the moment that $P_u\cap P_v$ contains a lattice point
$x$ in its relative interior. Then $x-v\in\L_{P_u}$ and
$x-u\in\L_{P_v}$. Assume in addition that $x-v$ is in the relative
interior of $P_u$. Then $x-v-u\in\L_P$. The automorphism
$[e_u^\lambda,e_v^\mu]$ fixes pointwise both $\L_{P_u}$ and
$\L_{P_v}$. On the other hand $\gp(S_P)$ is generated by
$\L_{P_u}$ and $x-v-u$ (or, similarly, by $\L_{P_v}$ and $x-v-u$).
Therefore, it suffices to show that the elementary automorphisms
$e_u^\lambda$ and $e_v^\mu$ commute at $x-v-u$. It is clear that
both automorphisms restrict to an automorphism of the subalgebra
$k[\Q]\subset k[P]$ where $\Q\subset P$ is the lattice 2
dimensional parallelepiped with vertices $x$, $x-v$, $x-u$,
$x-u-v$. (Clearly, $\Q$ is a unit lattice square up to unimodular
transformation.) Now the desired commutativity follows from the
easily checked equation
$$
(e_u^\lambda\circ e_v^\mu)\bigl(xu^{-1}v^{-1}\bigr)=
xu^{-1}v^{-1}+\lambda(xv^{-1})+\mu(xu^{-1})+\lambda\mu x
=(e_v^\mu\circ e_u^\lambda)(xu^{-1}v^{-1})
$$
where the semigroup operation is written multiplicatively.

Consider the general case. Since $S_P$ consists of
non-zero-divisors it is enough to show that $e_u^\lambda\circ
e_v^\mu$ and $e_v^\mu\circ e_u^\lambda$ coincide on Veronese
subalgebras $k[P]_{(c)}$ for all sufficiently big natural numbers
$c$. It is in general not true that $k[P]_{(c)}=k[cP]$ (unless $P$
is a normal polytope, that is the semigroup $S_P$ is normal, see
\cite{BrGT} for a detailed study of these properties). But the
monomials of degree $c$ in $k[P]_{(c)}$ constitute a subset of
$\L(cP)$ such that all the arguments above that we have used for
$P$ apply to it. It only remains to notice that for $c$ big enough
the existence of the desired lattice points in the relative
interiors as above is guaranteed.
\end{proof}

\section{The stable group of elementary automorphisms}\label{STABLE}

{\em From now on we assume that the polytopes being considered
have at least one column vector.}

For a polytope $P$ the group of elementary automorphisms
$\EE_R(P)$ may not be perfect. For instance, $P$ can have only one
column structure (like the example $P=\conv((0,0),(2,0),(1,2))$),
and then $\EE_R(P)$ is isomorphic to the additive group $R$.
Observe that such a polytope is automatically balanced.

We resolve this difficulty by the use of {\it doubling spectra} of
polytopes.

\begin{definition}\label{dspec}
An ascending infinite chain of polytopes
$\Pp=(P=P_0\subset P_1\subset\ldots)$ is called a {\it doubling
spectrum} if the following conditions hold:
\begin{itemize}
\item[(i)]
for every $i\in\ZZ_+$ there exists a column vector
$v\subset\Col(P_i)$ such that $P_{i+1}=P_i^{\sq_v}$,
\item[(ii)]
for all $i\in\ZZ_+$ and every $v\in\Col(P_i)$ there is an index
$j\geq i$ such that $P_{j+1}=P_j^{\sq_v}$.
\end{itemize}
\end{definition}

Here we have used Lemma \ref{colext} which (together with
condition (i)) allows one to consider $v$ as an element of
$\Col(P_j)$. The second part of the definition simplifies the
construction of doubling spectra. For example, if
$v,-v\in\Col(P)$, then it does not matter whether we double along
$P_v$ or $P_{-v}$ since there exist an automorphism of $P$
exchanging $v$ and $-v$.

We say that $v\in\Col(P_i)$ is {\it decomposed} at the $j$th step
in $\Pp$ for some $j\geq i$ if $P_{j+1}=P_j^{\sq_v}$. We will need
the fact that every column vector, showing up in a doubling
spectrum, gets decomposed infinitely many times:

\begin{lemma}\label{infdec}
For all $i\in\ZZ_+$ every vector $v\in\Col(P_i)$ is decomposed at
infinitely many steps in $\Pp$.
\end{lemma}

\begin{proof}
The desired infinite series of decompositions is derived as
follows. First fix some $j_1$ such that
$P_{j_1+1}=P_{j_1}^{\sq_v}$ and then choose $j_2<j_3<\cdots$
recursively such that
$$
P_{j_{t+1}+1}=\sqd\Bigl(P_{j_{t+1}},
  \delta^+(v,P_{j_t})\Bigr).
$$
The vector $v$ gets decomposed at each of the indices $j_t$
because $\delta^+(v,P_{j_t})$ and $v$ share the base facet in
$P_{j_{t+1}}$.
\end{proof}

Associated to a doubling spectrum $\Pp$ is the ``infinite
polytopal'' algebra
$$
R[\Pp]=\lim_{i\to\infty}R[P_i]
$$
and the filtered union
$$
\Col(\Pp)=\lim_{i\to\infty}\Col(P_i).
$$
The embeddings meant here are induced by the facet embeddings
$P_i\subset P_{i+1}$ which, by Lemma \ref{colext} imply natural
embeddings $\Col(P_i)\subset\Col(P_{i+1})$. For convenience of
notation we will often identify $P_i$ and $\Col(P_i)$ with their
images under these embeddings, without mentioning this explicitly.

All elements $v\in\Col(\Pp)$ and $\lambda\in R$ give rise to a
graded automorphism of $R[\Pp]$ as follows: we choose an index $i$
big enough so that $v\in\Col(P_i)$. Then the elementary
automorphisms $e_v^{\lambda}\in\EE_R(P_j)$, $j\geq i$ constitute a
compatible system and, therefore, define a graded automorphism of
$R[\Pp]$. This automorphism will again be called `elementary' and
it will also be denoted by $e_v^{\lambda}$.

Let $\EE(R,\Pp)$ denote the subgroup of the automorphism group of
$R[\Pp]$ that is generated by the elementary automorphisms. We
will call $\EE(R,\Pp)$ the {\it stable group of elementary
automorphisms} over $P$. Clearly, there are many doubling spectra
starting from $P$, but all the resulting stable groups are
pairwise naturally isomorphic.

\begin{proposition}\label{specinv}
Let $\Pp=(P\subset P_1\subset P_2\subset\cdots)$ and $\Qq=P\subset
Q_1\subset Q_2\subset\cdots)$ be two doubling spectra. Then the
groups $\EE(R,\Pp)$ and $\EE(R,\Qq)$ are naturally isomorphic.
\end{proposition}

\begin{proof}
Consider the `infinite lattice polytopes'
$$
{\cal P}=\bigcup_{\ZZ_+}P_i\subset\bigoplus_{\ZZ_+}\RR\ \
\text{and}\ \ {\cal
Q}=\bigcup_{\ZZ_+}Q_i\subset\bigoplus_{\ZZ_+}\RR,
$$
where we mean the filtered unions of polytopes and their ambient
Euclidean spaces. It is enough to show that there is a global
affine (i.~e.\ linear + constant) transformation
$$
\Theta:\bigoplus_{\ZZ_+}\RR\to\bigoplus_{\ZZ_+}\RR,\qquad
\Theta\bigl(\bigcup_{\ZZ_+}\L_{P_i}\bigr)=
\bigcup_{\ZZ_+}\L_{Q_i},
$$
such that $\Theta({\cal P})={\cal Q}$ and
$\Theta\bigl(\Col(\Pp)\bigr)=\Col(\Qq)$.

The existence of such $\Theta$ is established as follows. By Lemma
\ref{colbal} (and Definition \ref{dspec}) for every index
$i\in\ZZ_+$ there exists $j_i\geq i$ such that $\Col(P_i)$ is a
subset of $\Col(Q_{j_i})$. Consider the minimal subpolytope
$Q'\subset Q_{j_i}$ for which $P\subset Q'$ and
$\Col(P_i)\subset\Col(Q')$. It is clear from Proposition
\ref{welldef} and Definition \ref{dspec} that $P_i$ and $Q'$ are
naturally isomorphic as lattice polytopes and, moreover, $Q'$ is a
face of $Q_{j_i}$. In particular we obtain an isomorphism between
$P_i$ and a face of $Q_{j_i}$. Let $\Theta_i$ denote the affine
continuation of this isomorphism to the ambient Euclidean spaces.
Then the $\Theta_i$ constitute a compatible system and, therefore,
induce a global affine transformation
$\Theta:\bigoplus_{\ZZ_+}\RR\to\bigoplus_{\ZZ_+}\RR$. It is an
easy exercise to show that $\Theta$ satisfies all the desired
conditions.
\end{proof}

As a consequence of Proposition \ref{specinv} we can use the
notation $\EE(R,P)$ for $\EE(R,\Pp)$ where $\Pp$ is some doubling
spectrum starting with $P$. In the sequel we will assume that for
a polytope $P$ we have fixed a doubling spectrum $\frak P$.

\begin{remark}\label{isosequ}
One can define the group $\EE(R,P)$ using sequences of polytopes
$\Pp'=(P=P'_0\subset P_1'\subset \cdots)$ that are more general
than doubling spectra. In particular, suppose that
$\Pp=(P=P_0\subset P_1\subset\cdots)$ is a doubling spectrum for
$P$ and $\Pp'=(P_0'\subset P_1'\subset \cdots)$ is a sequence of
polytopes for which there exist isomorphisms $\phi_i:P_i\to P_i'$
of polytopes that commute with the embeddings $P_i\subset P_{i+1}$
and $P_i'\to P_{i+1}'$. Then $\Pp'$ need not be a doubling
spectrum in the strict sense since condition (ii) is not invariant
under isomorphisms as just described. However, there evidently
exists a natural isomorphism $\EE(R,\Pp)\approx \EE(R,\Pp')$.

For instance, if we start from the unimodular simplex $\Delta_n$,
$n\in\NN$, and consider the sequence $\Pp'=(\Delta_n=P'_0\subset
P_1'\subset\cdots)$, in which $P'_{i+1}={P'_i}^{\sq_v}$,
$i\in\NN$, for the same column vector $ v\in\Col(\Delta_n)$, then
the resulting sequence of unstable groups is naturally identified
with the familiar sequence of groups of elementary matrices
$$
E_{n+1}(R)\subset E_{n+2}(R)\subset\cdots, \qquad*\mapsto
\begin{pmatrix}
*&0\\
0&1
\end{pmatrix}.
$$
In particular, $\EE(R,\Delta_n)=\E(R)$ for all $n\in\NN$.
\end{remark}

Let $u,v\in\Col(\Pp)$. It is easily observed that the
$\la(P_i)_v,u\ra$ are the same for all $P_i$ from $\Pp$ such that
$u,v\in\Col(P_i)$. We will use the notation $\la\Pp_v,u\ra$ for
this common value.

\begin{proposition}\label{perf}
Let $P$ be a polytope, $u,v\in\Col(\Pp)$ and $\lambda,\mu\in R$.
Then:
\begin{itemize}
\item[(a)]
$e_u^{\lambda}\circ e_u^{\mu}=e_u^{\lambda+\mu}$.
\item[(b)]
The exact analogues of the equations in Theorem \ref{comrel} hold
once $\la P_v,u\ra$ is changed to $\la\Pp_v,u\ra$.
\item[(c)]
$\EE(R,P)$ is perfect.
\end{itemize}
\end{proposition}

(a) and (b) are proved by an easy reduction to $R[P_i]$ with $i$
big enough and (c) follows from (b) by the equation
$(\ref{DOUBL}_2)$. In particular, the perfectness does not depend
on whether or not $P$ is balanced.

The following remark shows that naive analogies between the stable
group of elementary matrices $\E(R)$ and $\EE(R,P)$ may fail.

\begin{remark}\label{nfiltun}
(a) In the special case of unit simplices we get the stable group
of elementary matrices over $R$: $\EE(R,\Delta_n)=\E(R)$ for all
$n\in\NN$.

(b) In general the groups $\EE_R(P)$ are not perfect. However,
after finitely many steps in the doubling spectrum one arrives at
a polytope $P''$ for which this group is perfect.

In fact, after finitely many steps each base facet in $P$ has been
used for a doubling, and in the polytope $P'$ then constructed
each base facet has an invertible column vector. This property is
preserved under further doublings. Therefore all the vectors
$\delta^+$ and $\delta^-$ that come up in doublings of $P'$ are
automatically decomposed, and after finitely many doublings
starting from $P'$ one arrives at a polytope $P''$ in which all
the column vectors of $P'$, and thus all of $P''$, are decomposed.

(c) The group $\EE(R,P)$ is in general {\it not} the filtered
union of the unstable groups $\EE_R(P_i)$. Consider the simple
example of the segment $2\Delta_1$. Then the second term in $\Pp$
can be identified with the triangle
$2\Delta_2=\conv\bigl((0,0),(2,0),(0,2)\bigr)\subset\RR^2$ so that
$2\Delta_1$ is the lower edge (see Figure \ref{NonEmb}).
\begin{figure}[htb]
\begin{center}
\psset{unit=1cm}
\def\vertex{\pscircle[fillstyle=solid,fillcolor=black]{0.07}}
\begin{pspicture}(1.0,0)(1,2.2)
 \footnotesize
 \pspolygon[style=fyp](0,0)(2,0)(0,2)
 \multirput(0,0)(0,1){3}{\vertex}
 \multirput(1,0)(0,1){2}{\vertex}
 \multirput(2,0)(0,1){1}{\vertex}
 \psline[linewidth=1pt]{->}(1,0)(0,0)
 \psline[linewidth=1pt]{->}(1,0)(2,0)
 \rput(0.5,-0.3){$-v$}
 \rput(1.5,-0.3){$v$}
\end{pspicture}
\caption{The polytope $2\Delta_2$} \label{NonEmb}
\end{center}
\end{figure}
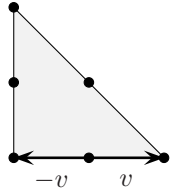
Consider the vectors  $v=(1,0)$ and $-v=(-1,0)$ from
$\Col(2\Delta_1)$. Assume $2\neq 0$ in $R$. Then the element
$\epsilon=(e_v^1\circ e^{-1}_{-v}\circ e_v^1)^2\in\EE(R,\Delta_2)$
is not the identity automorphism of $R[\Pp]$ (it switches signs on
the second layer of $2\Delta_2$) whereas the restriction of
$\epsilon$ to $R[2\Delta_1]$ is the identity automorphism. In
particular there is no natural group embedding
$\EE_R(2\Delta_1)\to\EE(R,2\Delta_1)$.

(d) On the other hand, every element $e\in\EE_R(P_i)$, $i\in\ZZ_+$
is a restriction to $R[P_i]$ of some element of $\EE(R,P)$ and
every element $\epsilon\in\EE(R,P)$ restricts to an element of
$E_R(P_i)$ whenever $i$ is big enough. Clearly, we have the
following approximation principle: if two elements
$\epsilon,\epsilon'\in\EE(R,P)$ restrict to the same elements of
$\EE_R(P_i)$ for all sufficiently large $i$ then
$\epsilon=\epsilon'$.

(e) Unlike the group $\EE(R,P)$, the Steinberg group $\St(R,P)$,
introduced in Section \ref{SCHUR}, is the direct limit of the
corresponding unstable groups (see Remark \ref{filtun} below).
\end{remark}

The next results is just a standard fact in the classical
situation of elementary matrices. However, the proof is no longer
straightforward in the polytopal case, not even for balanced
polytopes.

\begin{theorem}\label{trivcen}
Let $R$ be a ring and $P$ a polytope (not necessarily balanced).
Then the center $Z\bigl(\EE(R,P)\bigr)$ is trivial.
\end{theorem}

\begin{proof}
Choose $\epsilon\in Z\bigl(\EE(R,P)\bigr)$. There is no loss of
generality in assuming $\epsilon=e_{v_1}^{\lambda_1}\circ
\cdots\circ e_{v_k}^{\lambda_k}$ for some $k\geq0$,
$\lambda_1,\ldots,\lambda_k\in R$ and $v_1,\ldots,v_k\in\Col(P)$.
In particular, $\epsilon$ restricts to an automorphism of each of
the $R[P_i]$, $i\in\ZZ_+$.
\smallskip

\noindent {\it Step 1.} We claim that for every $x\in \L_P$ there
exists $a_x\in R^*$ such that $\epsilon(x)=a_xx$. ($R^*$ denotes
to the multiplicative group of invertible elements.)

Let $F\subset P$ be a facet. It is enough to show
that every element $x\in\L_P$ satisfies the condition
$\epsilon(x)=a_1x_1+\cdots+a_sx_s$ for some $a_i\in R$ and
$x_i\in\L_P$ with $\het_F(x_i)\geq\het_F(x)$. In fact, by running
$F$ through $\FF(P)$ we prove the claim.

In the argument below we consider the polytopes $P^{\sq_F}$ for
all facets $F\subset P$. Clearly, if $F$ is a not a base facet of
an element of $\Col(P)$, then $P^{\sq_F}$ does not show up in a
doubling spectrum starting with $P$ -- we use this doubled
polytope only as an auxiliary object.

If $F=P_v$ for some $v\in\Col(P)$, then the two automorphisms
$\epsilon,e_{\delta^-(F,P)}^1\in\EE(R,P)$ commute. If $F$ is not
the base facet of any of the column vectors, then none of the
products $\delta_{F,P}^-v$ and $v\delta_{F,P}^-$ exists for any
$v\in\Col(P)$. In fact, we have $\la
P^{\sq_F}_v,\delta_{F,P}^-\ra=0$ (by equation $(\ref{DOUBL}_1)$)
and $\la P_{\delta^-(F,P)},v\ra= \la P^-,v\ra=0$. Since
$v+\delta_{F,P}^-\neq 0$, Theorem \ref{comrel} implies that
$\epsilon$ and $e_{\delta^-(F,P)}^1$ again commute. For simplicity
put $\delta^-=\delta_{F,P}^-$.

Consider the point $x^\vt \in P^\vt $ that corresponds to $x$. We
have
$$
\epsilon(x^\vt )\in R\bigl[P^{\sq_F}\bigr]=b_1y_1+\cdots+b_ty_t
$$
for some $b_j\in R$ and $y_j\in\L_{P^{\sq_F}}$ with
$\het_{P^-}(y_j) =\het_F(x)$, where $\het_F(x)$ is considered in
$P$. It is easily seen that
$$
(e^1_{\delta^-}\circ\epsilon)(x^\vt
)=b_1z_1+\cdots+b_tz_t+(\text{$R$-linear combination of points $u$
with}\ \het_{P^-}(u)>0),
$$
where $z_j=y_j+\het_F(x)\cdot\delta^-\in\L_{P^-}$. Clearly, $\het_F(z_j)\geq
\het_F(x)$.

On the other hand
$$
(\epsilon\circ e^1_{\delta^-})(x^\vt
)=\epsilon(x)+(\text{$R$-linear combination of points $u$ with}\
\het_{P^-}(u)>0)
$$
Since $\epsilon$ and $\delta^-$ commute, we get what we want by
comparing these two \smallskip expressions.

\noindent{\it Step 2.} Let $Q$ be a polytope,
$w\in\Col(Q)$ and $\alpha$ be an $R$-automorphism of $R[Q]$ of
type $y\mapsto a_yy$, $a_y\in R^*$ for all $y\in\L_Q$. Put
$$
\alpha(w)=\frac{\alpha(z\cdot w)}{\alpha(z)}
$$
where $z\in\L_Q\setminus\L_{Q_w}$ and the operation in $\ZZ^n$ is
written multiplicatively. The ratio is in fact independent of the
choice of $z$ because $\alpha(z'w)\alpha(z)=\alpha(zw)\alpha(z')$.
Then
\begin{equation}
\tag{$*$}\alpha\circ e_w^1\circ\alpha^{-1}=e_w^{\alpha(w)}.
\end{equation}
(See also Lemma 4.4(b) in [BG2].)

Returning to our central automorphism $\epsilon$ we see that
$a_x=a_{x'}$ whenever the points $x,x'\in\L_P$ differ by a vector
from $\Col(P)$. Since we can replace $P$ by some polytope obtained
by successive doublings at base facets of column vectors, it is
enough to apply the following lemma: clearly $a_y=1$ for the
lattice point $y$ constructed there, and thus $a_x=1$ since $x$
and $y$ differ by a sum of column vectors.
\end{proof}

\begin{lemma}\label{facetize}
Let $P$ be a polytope, $x\in \L_P$, and $v_1,\dots,v_k$
be column vectors of $P$. Set $P_0=P$, $P_1=P^{\sq_{v_1}}$,
$P_2=(P_1)^{\sq_{v_2}},\dots$. Then there exists $y\in \L_Q$,
$Q=P_k$, such that $x$ and $y$ differ by a sum of column vectors
of $Q$, and $y\in Q_{v_i}$ for $i=1,\dots,k$.
\end{lemma}

\begin{proof}
We construct a chain $x=x_0,\dots,x_k=y$ of lattice points $x_i\in
P_i$ such that $x_i-x_{i-1}$ is a multiple of $\delta_{v_i}^+$ and
$x_i\in (P_i)_{v_i}$ for $i=0,\dots,k-1$.

Let $i\ge1$. If $x_{i-1}\in (P_{i-1})_{v_i}$, then we choose
$x_i=x_{i-1}$. Otherwise $(P_{i-1})_{v_i} \neq(P_{i-1})_{v_j}$ for
$j=1,\dots,i-1$, and we set
$x_i=x_{i-1}+\het_{v_i}(x_{i-1})\delta_{v_i}^+$. Since
$\delta_{v_i}^+$ is parallel to the extension of $(P_{i-1})_{v_j}$
to $P_i$ for $j=1,\dots,i-1$ (by $(\ref{DOUBL}_1)$) and
$´\het_{v_i}\bigl(x_{i-1}+\het_{v_i}(x_{i-1})\delta_{v_i}^+\bigr)=0$,
we have reached a lattice point with the desired properties.
\end{proof}

\section{Schur multiplier of $\EE(R,P)$}\label{SCHUR}

Next we proceed in analogy with the ordinary algebraic $K$-theory and
define the {\it polytopal Steinberg groups}. In the proofs below we use
a number of modifications of Milnor's arguments \cite[\S 5]{M}. Many of
the difficulties that show up for general balanced polytopes are invisible
in the special case of unit simplices.

Fix a doubling spectrum $\Pp$ starting with $P$. The {\it
unstable} Steinberg groups $\St_R(P_i)$ are defined as the groups
generated by the $x_v^\lambda$, $v\in\Col(P_i)$, $\lambda\in R$,
which are subject to the relations
$$
x_v^\lambda x_v^\mu=x_v^{\lambda+\mu}
$$
and
$$
[x_u^\lambda,x_v^\mu]=
\begin{cases}
\Pi_{i=1}^n x_{u+iv}^{-\binom{n}{i}\lambda\mu^i}&
\text{if}\ uv\ \text{exists, } n=\la P_v,u\ra,\\
\\
\Pi_{i=1}^n x_{v+iu}^{\binom{n}{i}\mu\lambda^i}&
\text{if}\ vu\ \text{exists, } n=\la P_v,u\ra,\\
\\
1&\text{if}\ u+v\notin\Col(P)\cup\{0\}.
\end{cases}
$$
The {\it stable Steinberg group} $\St(R,P)$ is defined through the
obvious analogues of the equations above, where $\la P_v,u\ra$ and
$\la P_u,v\ra$ are correspondingly changed to $\la\Pp_v,u\ra$ and
$\la\Pp_u,v\ra$. Clearly, arguments similar to those in
Proposition \ref{specinv} show that $\St(R,P)$ does not depend on
the choice of the doubling spectrum $\Pp$.

\begin{remark}\label{filtun}
(a) There are natural surjective homomorphisms
$\St_R(P_i)\to\EE_R(P_i)$ and $\St(R,P)\to\EE(R,P)$. Moreover,
contrary to the groups of elementary automorphisms (see Remark
\ref{nfiltun}), one has natural group homomorphisms between
successive unstable Steinberg groups, and
$$
\St(R,P)=\lim_{\to}\St_R(P_i).
$$
However, these homomorphisms may be non-injective, even in the
classical situation $P=\Delta_n$ -- here we enter the topic of
injective $K_2$-stabilization.

(b) The group $\St(R,P)$ is always perfect, like $\EE(R,P)$.

(c) Remark \ref{isosequ} applies here as well: $\St(R,P)$ can be
computed from any ascending sequence of lattice polytopes that is
isomorphic to a doubling spectrum.
\end{remark}

\begin{proposition}\label{centker}
Let $P$ be a balanced polytope and $R$ be a ring. Then
then the center $Z(\St(R,P))$ is the kernel of the natural
homomorphism $\St(R,P)\to\EE(R,P)$.
\end{proposition}

Later on $\Ker(\St(R,P)\to\EE(R,P))$ for $P$ balanced will be called the
{\it polytopal Milnor group} of $R$ corresponding to $P$. We denote it
by $K_2(R,P)$.

\begin{proof}[Proof of Proposition \ref{centker}]
As remarked, we have fixed a doubling spectrum $\Pp=(P\subset
P_1\subset\cdots)$. For every $i\in\NN$ consider the subsets
\begin{align*}
U^{i+1}&=\{u\in\Col(P_{i+1})\ |\ \la P_i,u\ra=1\},\\
V^{i+1}&=\{v\in\Col(P_{i+1})\ |\ \la P_i,v\ra=-1\}
\end{align*}
of $\Col(P_{i+1})$. (Unlike the classical situation the two sets
$U^{i+1}$ and $V^{j+1}$ are not completely similar.)

We have the subgroups ${\frak U}^{i+1}, {\frak
V}^{i+1}\subset\St_R(P_{i+1})$ generated correspondingly by the
$x_u^\lambda$, $v\in U^{i+1}$ and $x_v^\mu$, $v\in V^{i+1}$
($\lambda,\mu\in R$). Assume $U^{i+1}=\{u_1,\cdots,u_r\},\ \
r=\#U^{i+1}$ and $V^{i+1}=\{v_1,\cdots,v_s\},\ \ s=\#V^{i+1}$.
\medskip

\noindent{\bf Claim.} The mappings
\begin{alignat*}{2}
R^r&\to{\frak U}^{i+1},&\qquad(\lambda_1,\cdots,\lambda_r)&\mapsto
x_{u_1}^{\lambda_1}\cdots x_{u_r}^{\lambda_r},\\
R^s&\to{\frak V}^{i+1},&\qquad (\mu_1,\cdots,\mu_s)&\mapsto
x_{v_1}^{\mu_1}\cdots x_{v_s}^{\mu_s}
\end{alignat*}
are  group isomorphisms, where the free $R$-modules $R^r$ and
$R^s$ are viewed as additive abelian groups. Moreover, the
canonical surjection $\pi:\St_R(P_{i+1})\to\EE_R(P_{i+1})$ is
injective both on ${\frak U}^{i+1}$ and ${\frak V}^{i+1}$.
\medskip

First one observes that $ww'$ does not exist for any pair
$$
(w,w')\in\bigl(U^{i+1}\times U^{i+1}\bigr)\cup{\big
(}V^{i+1}\times V^{i+1}\bigr).
$$
Lemma \ref{colv} implies this for the pairs from $V^{i+1}\times
V^{i+1}$ without the condition that $P$ is balanced, whereas this
condition is needed for the pairs from $U^{i+1}\times U^{i+1}$: if
$ww'$ existed then for some $(w,w')\in U^{i+1}\times U^{i+1}$ then
$\la P_i,ww'\ra=\la(P_{i+1})_{\delta^-},ww'\ra=2$, contradicting the
condition that $P$ is balanced.

As a result of this observation all elements $\sigma\in{\frak
U}^{i+1}$ and $\rho\in{\frak V}^{i+1}$ have presentations
$$
\text{(1)}\quad\sigma=x_{u_1}^{\lambda_1}\cdots
x_{u_r}^{\lambda_r}\qquad\text{and}\qquad
\text{(2)}\quad\rho=x_{v_1}^{\mu_1}\cdots x_{v_s}^{\mu_s}.
$$
For all points $x\in\L_{P_i}$ and $y\in\L_{P_{i+1}}$ with
$\het_{P_i}(y)=1$ we have
\begin{align*}
\bigl(\pi(\sigma)\bigr)(x)&=x+\tau_1\lambda_1(xu_1)+\cdots+
\tau_r\lambda_r(xu_r)+f_{\geq2}\\
\bigl(\pi(\rho)\bigr)(y)&=y+\mu_1(yv_1)+\cdots+\mu_s(yv_s)
\end{align*}
where $f_{\geq2}$ is an $R$-linear combination of the points from
$L_{P_{i+1}}$ with height above $P_i$ at least 2 and
$\tau_i\in\{0,1\}$ for $i\in[1,r]$. (The semigroup operation is
written multiplicatively). It is immediate from the second
equation that the presentation (2) is uniquely determined by
$\pi(\rho)$. Observe that for each $i\in[1,r]$ there exists
$x\in\L_{P_i}$ such that the corresponding $\tau_i$ is 1.
Therefore, by running $x$ through $\L_{P_i}$ we see that likewise
$\pi(\sigma)$ uniquely determines the presentation (1). Now the
claim follows.

Since $P$ is balanced, for all vectors $u\in U^{i+1}$, $v\in V^{i+1}$
and $w\in\Col(P_i)$ and elements $\lambda,\mu,\nu\in F$ we have:
$$
x_w^\nu x_u^\lambda x_w^{-\nu}=
\begin{cases}x_{wu}^{-\nu\lambda}x_u^\lambda&\text{if}\ wu\ \text{exists},\\
x_{uw}^{\lambda\nu}x_u^{\lambda}&\text{if}\ uw\ \text{exists},\\
x_u^{\lambda}&\text{if}\ wu\notin\Col(P),
\end{cases}
$$
and
$$
x_w^\nu x_v^{\mu}x_w^{-\nu}=
\begin{cases}x_{vw}^{\mu\nu}x_v^\mu&\text{if}\ vw\ \text{exists},\\
x_v^\mu&\text{if}\ vw\ \text{does not exist}.
\end{cases}
$$
The listed cases exhaust all possibilities since $u+w\neq 0$,
$v+w\neq0$ and the product $wv$ does not exist -- we have
$\la(P_{i+1})_v,w\ra=\la P_i,w\ra=0$ and Proposition
\ref{maincrit}(a) applies. (It is not difficult to find examples
for which $uw$ exists.)

Since $wu\in U^{i+1}$, $uw\in U^{i+1}$ and $vw\in V^{i+1}$
whenever the corresponding product exists, we arrive at the
inclusions:
$$
z{\frak U}^{i+1}z^{-1}\subset{\frak U}^{i+1}\quad\text{and}\quad
z{\frak V}^{i+1}z^{-1}\subset{\frak V}^{i+1}.
$$
for all $z\in\St_R(P_i)$. Changing $z$ by $z^{-1}$ we get
\begin{equation}
\tag{3} z{\frak U}^{i+1}z^{-1}={\frak U}^{i+1}\quad\text{and}\quad
z{\frak V}^{i+1}z^{-1}={\frak V}^{i+1}.
\end{equation}
for all $z\in\St_R(P_i)$.

Now we are ready to prove that
$$
\Ker\bigl(\pi:\St(R,P)\to\EE(R,P))= Z(\St(R,P)\bigr).
$$
Choose an element $z\in\St(R,P)$ is such that $\pi(z)=1$. We want
to show that $z x_v^\lambda z^{-1}=x_v^\lambda$ for all
$v\in\Col(\Pp)$ and $\lambda\in R$. Since $\St(R,P)$ is the
inductive limit of the groups $\St_R(P_i)$ and $\pi(z)$ restricts
to an automorphism of $R[P_i]$ for $i\gg0$, there exist
$z_i\in\Ker\bigl(\St_R(P_i)\to\E_R(P_i)\bigr)$ mapping to $z$
provided $i\gg0$ and such that $z_{i+1}$ is the image of $z_i$ in
$\St_R(P_{i+1})$. We can also assume $v\in\Col(P_i)$ and that $v$
gets decomposed in $\Pp$ at $P_i$. We will show that $z_i$ maps to
the center
$Z\bigl(\Im\bigl(\St_R(P_i)\to\St_R(P_{i+1})\bigr)\bigr)$.
Clearly, this implies that $z$ is central in $\St(R,P)$.

For the elements of $\St_R(P_{i+1})$ we have
$$
x_v^\lambda=\bigl[x^1_{\delta^+(v,P_i)},x_{v^\vt
}^\lambda\bigr],\quad \ x^1_{\delta^+(v,P_i)}\in{\frak
U}^{i+1},\quad x_{v^\vt }^\lambda \in{\frak V}^{i+1},
$$
where $v^\vt \in\Col(P_{i+1})$ is the column vector corresponding
to $v$. The claim above and (3) imply that $z_i$ commutes both with
elements of ${\frak U}^{i+1}$ and ${\frak V}^{i+1}$. Hence it also
commutes with $x_v^{\lambda}$.

We have shown the inclusion
$\Ker\bigl(\pi:\St(R,P)\to\EE(R,P)\bigr) \subset Z(\St(R,P))$. The
other inclusion follows from Theorem \ref{trivcen}.
\end{proof}

Before the discussion of balanced polytopes in general let us
treat a polytope that represents the case in which there are three
column vectors with exactly one product.

\begin{lemma}\label{specquad}
Let $P=\conv\bigl((0,0),(3,0),(1,2),(0,1)\bigr)$ and $R$ be a
ring. Then the group $\St(R,P)$ is a universal central extension of
$\EE(R,P)$.
\end{lemma}

\begin{proof}
We have $\Col(P)=\{u,v,w\}$ for $u=(0,-1)$, $v=(1,0)$ and
$w=(1,-1)$. The product table for $\Col(P)$ consists of only the
equation $uv=w$ (see Figure \ref{OneProd}).
\begin{figure}[htb]
\begin{center}
 \psset{unit=1cm}
 \def\vertex{\pscircle[fillstyle=solid,fillcolor=black]{0.05}}
\begin{pspicture}(-0.3,0)(4,2.5)
 \pspolygon[style=fyp,linecolor=darkgray](0,0)(3,0)(1,2)(0,1)
 \footnotesize
 \multirput(0,1)(1,0){3}{\vertex}
 \multirput(0,0)(1,0){4}{\vertex}
 \rput(1,2){\vertex}
 \psline[style=fatline]{->}(0,1)(1,0)
 \psline[style=fatline]{->}(0,1)(0,0)
 \psline[style=fatline]{->}(0,0)(1,0)
 \rput(-0.2,0.6){$u$}
 \rput(0.6,-0.2){$v$}
 \rput(0.7,0.7){$w$}
\end{pspicture}
\end{center}
\caption{Three column vectors with one product} \label{OneProd}
\end{figure}
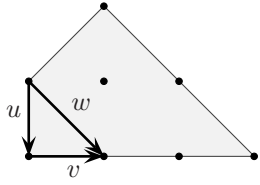
We now describe the group $\St(R,P)$ explicitly.

First we make the following observation: if $Q$ is a polytope and
$v,-v\in\Col(Q)$ then $Q^{\sq_v}$ is naturally isomorphic to
$Q^{\sq_{-v}}$ so that the points of $Q$ are mapped to themselves
under this isomorphism. It follows from this observation and Lemma
\ref{colbal} that for producing a sequence of polytopes isomorphic
to a doubling spectrum (starting with $Q$) one only needs to
decompose at each step one of the original column vectors of $Q$
or one of the vectors of type either $y^\vt $ or $\delta^+$. (See
Remarks \ref{isosequ} and \ref{filtun}.)

Fix a doubling spectrum $\Pp=(P\subset P_{1}\subset\cdots)$.
Consider the points $a=(0,1)$, $b=(0,0)$ and $c=(1,0)$ from
$\L_{P}$. We have $a+u=b$ and $b+v=c$. Based on the aforementioned
general observation one easily sees that for every $i$ the set
$\Col(P_i)$ looks as follows. There is certain system of points
in $\L_{P_i}$
$$
a=x_1,x_2\ldots,x_{j_0}=b,x_{j_0+1},\ldots,x_m=c,\qquad 1<j_0<m
$$
such that
$$
x_{j'}-x_{j''}\in\Col(P_i)\iff\ (j',j'')\in \bigl(Z_1\cup Z_2\cup
Z_3\cup Z_4\bigr)\setminus\{(1,1),\ldots,(m,m)\},
$$
where
\begin{alignat*}{2}
Z_1&=[1,j_0-1]\times[1,j_0-1]\qquad&
Z_2&=[j_0,m]\times[1,j_0-1],\\
Z_3&=[j_0,m-1]\times[j_0,m-1]\qquad& Z_4&=\{m\}\times[1,m-1].
\end{alignat*}
Moreover, the numbers $m$, $j_0$ and $m-j_0$ can be arbitrarily
big if $i$ is big enough.

Now $\St(R,P)$ admits the following description. Let $A$ and $B$
denote two disjoint copies of $\NN$ and let $0$ be an ``origin'',
$0\notin A\cup B$. Then $\St(R,P)$ is generated by symbols
$x_{ij}^\lambda$ where $\lambda\in R$, $i\neq j$ and
$$
(i,j)\in\bigl(A\times A\bigr)\cup\bigl(B\times A\bigr)\cup\bigl(
B\times B\bigr)\cup\bigl(\{0\}\times A\bigr)\cup\bigl(\{0\}\times
B\bigr),
$$
and these symbols are subject to the standard Steinberg relations:
$$
x_{ij}^{\lambda}x_{ij}^{\mu}=x_{ij}^{\lambda+\mu}
\qquad\text{and}\qquad [x_{ij}^\lambda,x_{kl}^\mu]=
\begin{cases}
x^{\lambda\mu}_{il}&\text{if}\ j=k\ \text{and}\ i\neq l,\\
1&\text{if}\ j\neq k\ \text{and}\ i\neq l.
\end{cases}
$$
In Figure \ref{GenSt} we have tried to visualize these data.
\begin{figure}[htb]
\begin{center}
 \psset{unit=1cm}
 \def\vertex{\pscircle[fillstyle=solid,fillcolor=black]{0.05}}
\begin{pspicture}(-3,-0.5)(3,3)
 \pscircle[fillstyle=solid,fillcolor=light](-2,2){1}
 \pscircle[fillstyle=solid,fillcolor=light](2,2){1}
 \rput(0,-0.5){\vertex}
 \rput(-2.5,1.57){\vertex}
 \rput(-1.5,1.57){\vertex}
 \rput(-2,2.43){\vertex}
 \rput(2.5,1.57){\vertex}
 \rput(1.5,1.57){\vertex}
 \rput(2,2.43){\vertex}
 \psline[style=fatline]{->}(-1.5,1.57)(0,-0.5)
 \psline[style=fatline]{->}(1.5,1.57)(0,-0.5)
 \psline[style=fatline](1.5,1.57)(-1.5,1.57)
 \psline[style=fatline]{->}(1.5,1.57)(0,1.57)
 \multirput(-2.5,1.57)(4,0){2}{%
   \psline[style=fatline]{->}(0,0)(1,0)
   \psline[style=fatline]{->}(1,0)(0,0)
   \psline[style=fatline]{->}(0,0)(0.5,0.86)
   \psline[style=fatline]{->}(0.5,0.86)(0,0)
   \psline[style=fatline]{->}(1,0)(0.5,0.86)
   \psline[style=fatline]{->}(0.5,0.86)(1,0)
    }
\end{pspicture}
\end{center}
\caption{The generators of $\St(R,P)$} \label{GenSt}
\end{figure}
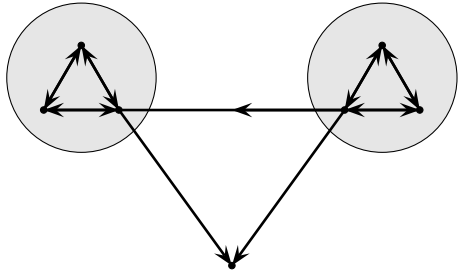

It only remains to notice that the proof of Theorem 5.10 in
\cite{M} goes through for $\St(R,P)$ without any change. It is of
course important that for a pair $(i,k)$ from our index set we can
always find $j$ such that $(i,j)$ and $(j,k)$ are also in the
index set.
\end{proof}

For general balanced polytopes $P$  one cannot visualize
$\St(R,P)$ in a similar way; see Example \ref{nomatrix} below.

\begin{theorem}\label{ucext}
Let $P$ be a balanced polytope and $R$ be a ring. Then
$\St(R,P)$ is a universal central extension of $\EE(R,P)$ and
$$
\Ker(\St(R,P)\to\EE(R,P))=Z(\St(R,P)).
$$
\end{theorem}

\begin{proof}
The second assertion is Proposition \ref{centker}.

Suppose
$$
\begin{CD}
1@>>>C@>>>Y@>\psi>>\St(R,P)@>>>1
\end{CD}
$$
is a central extension, i.~e. $C\subset Z(Y)$. We want to
construct a splitting homomorphism $\xi:\St(R,P)\to Y$. According
to \cite[\S5]{M} this will show the universality of the central
extension $\pi:\St(R,P)\to\EE(R,P)$.

Pick $x,x'\in\St(R,P)$. Since the extension $\psi$ is central, the
commutators of type $[y,y']$ with $y\in\psi^{-1}(x)$ and
$y'\in\psi^{-1}(x')$ coincide. Therefore we can use the notation
$[\psi^{-1}(x),\psi^{-1}(x')]$ for this common value.

For every generator $x_w^\lambda$ of $\St(R,P)$ with
$w\in\Col(P_i)$ and $\lambda\in R$ we choose an index $j\geq i$
such that the column vector $w$ is decomposed at $P_j$ and put
\begin{equation}
\tag{4} \xi(x_w^\lambda)=\bigl[\psi^{-1}(x^1_{\delta^+(w,P_j)}),
\psi^{-1}(x^\lambda_{w^\vt })\bigr].
\end{equation}
(Here $w^\vt \in\Col(P_{j+1})$ is the corresponding vector.) All
we have to show is that
\begin{itemize}
\item[(i)]
the $\xi(x_w^\lambda)$ are independent of the choices of $j$,
\item[(ii)]
the $\xi(x_w^\lambda)$ are subject to the same relations as
the $x_w^\lambda$.
\end{itemize}
To this end we make the following observation. Let $Q$ be the
quadrangle considered in Lemma \ref{specquad}, and
$u_0,v_0,w_0\in\Col(Q)$ be its column vectors. For some $k\ge 0$
Consider $u,v\in\Col(P_k)$ such that $uv$ exists. Then the
assignments $x^*_{u_0}\mapsto x^*_u$, $x^*_{v_0}\mapsto x^*_v$,
$x^*_{w_0}\mapsto x^*_{uv}$ ($*\in R$) give rise to a group
homomorphism $\St(R,Q)\to\St(R,P)$. This observation follows from
the definition of polytopal Steinberg groups and Lemma
\ref{colbal}: after the identifications $u\leftrightarrow u_0$,
$v\leftrightarrow v_0$ and $uv\leftrightarrow w_0$ the higher
members of the doubling spectrum $\Qq$ of $Q$ only admit column
vectors that also show up in $\Pp$. It is also crucial that the
defining relations for $\St(R,Q)$ are preserved by the
corresponding elements of $\St(R,P)$ because both $P$ and $Q$ are
balanced.

By Lemma \ref{specquad} every central extension of $\St(R,Q)$
splits. On the other hand the homomorphism $\Psi$ in the pull back
diagram
$$
\begin{CD}
Y\times_{\St(R,P)}\St(R,Q)@>\Psi>>\St(R,Q)\\
@VVV@VVV\\
Y@>>\psi>\St(R,P)
\end{CD}
$$
is central since $\psi$ is central. A splitting of $\Psi$ yields
the commutative triangle
$$
\begin{diagram}
    &         & \St(R,Q)&        & \\
    & \ldTo_f &           &  \rdTo & \\
Y   &         & \rTo^\psi &        & \St(R,P)\\
\end{diagram}
$$
which simultaneously shows several things:
\begin{itemize}
\item[(a)]
the definition of $\xi(x_w^\lambda)$ is in fact independent of the
choice of $j$,
\item[(b)]
$\xi(x_w^\lambda)\xi(x_w^\mu)= \xi(x_w^{\lambda+\mu})$ for all
$\lambda,\mu\in R$,
\item[(c)]
for all $\lambda,\mu\in R$ and all $w_1,w_2\in\Col(\Pp)$ for which
$w_1w_2$ exists one has
$[\xi(x_{w_1}^\lambda),\xi(x_{w_2}^\mu)]=\xi(x_{w_1w_2}^{\lambda\mu})$.
\end{itemize}
In fact, by lifting $\delta^+_{uv,P_k}$ and $(uv)^\vt $ to the
appropriate vectors in $\Col(\Pp_0)$ one observes that
\begin{equation}
\tag{5}
f(x_{w_0}^\lambda)={\big[}\psi^{-1}(x_u^1),\psi^{-1}(x_v^\lambda){\big]}
={\big [}\psi^{-1}(x^1_{\delta^+(uv,P_k)}),
\psi^{-1}(x^\lambda_{(uv)^\vt }){\big ]}
\end{equation}
for all $\lambda\in R$. Applying this formula to the column
vectors $w$, $\delta^+(w,P_j)$, and $w^\vt $ (playing the roles of
$uv$, $u$ and $v$) we prove the independence of $\xi(x_w^\lambda)$
of the choice of $j$. The other two properties of the
$\xi(x_v^\lambda)$ follow similarly by using appropriate mappings
$\St(R,Q)\to\St(R,P)$: for (b) we use the same mapping
$x_{u_0}^*\mapsto x_{\delta^+}^*$, $x_{v_0}\mapsto x_{w^\vt}^*$,
$x_{w_0}\mapsto x_w^*$
and for (c) we use the mapping
$x_{u_0}^*\mapsto x_{w_1}^*$, $x_{v_0}^*\mapsto x_{w_2}^*$,
$x_{w_0}^*\mapsto x_{w_1w_2}^*$, $*\in R$.

Only one Steinberg relation remains to be checked. Consider a pair
of column vectors $v,w\in\Col(P_i)$ such that
$v+w\notin\Col(P_i)\cup\{0\}$. We want to show that
$[\psi^{-1}(x_v^\lambda),\psi^{-1}(x_w^\mu)]=1$ for all
$\lambda,\mu\in R$. Like in the proof of Theorem \ref{comrel} one
must consider relations of this type separately.

Let $j$ be as in (4). Then
$$
{\big [}\psi^{-1}(x^1_{\delta^+(v,P_j)}),\psi^{-1}
(x^\lambda_{v^\vt }){\big ]}\in\psi^{-1}(x_v^\lambda).
$$
Assume we have shown that none of the products
\begin{equation}
\tag{5} \delta^+_{v,P_j}w,\quad w\delta^+_{v,P_j},\quad wv^\vt
,\quad\text{and}\quad v^\vt w
\end{equation}
exists. Then, since the corresponding sums are non-zero vectors,
we get
$$
\bigl[\psi^{-1}(x^1_{\delta^+(v,P_j)}),\psi^{-1}(x_w^\mu)\bigr],\
\bigl[\psi^{-1}(x^\lambda_{v^\vt }),\psi^{-1}(x_w^\mu)\bigr]\in
C\subset Z(Y).
$$
It follows that
$$
1=\Bigl[ \bigl[\psi^{-1}(x^1_{\delta^+(v,P_j)}),
\psi^{-1}(x^\lambda_{v^\vt })\bigr] ,\psi^{-1}(x_w^\mu)\Big]
=\bigl[\psi^{-1}(x_v^\lambda), \psi^{-1}(x_w^\mu)\bigr].
$$
It remains to show that the products (5) do not exist. (Notice
that this difficulty is absent in the case $P=\Delta_n$.)

For simplicity put $\delta^+=\delta^+_{v,P_j}$. Thus $\delta^+$
and $v$ have the same base facet in $P_{j+1}$, namely $P_j^\vt $.
(The arguments below use Proposition \ref{maincrit}(a) several
times.)

First observe that the inequality $\la P_j^\vt ,w\ra>0$ is
impossible because otherwise $wv$ would exist. This already
excludes the existence of $w\delta^+$. Also, the product $wv^\vt $
does not exist because $\la(P_{j+1})_{v^\vt },w\ra=\la
P_j^-,w\ra=0$.

By the equations $(\ref{DOUBL}_1)$ and $(\ref{DOUBL}_2)$ we have
$\la F,\delta^+\ra\leq0$ for every facet $F\subset P_{j+1}$
different from $P_j^-$. Since $w\in P_j^-$ the base facet
$(P_{j+1})_w\subset P_{j+1}$ is different from $P_j^-$. Therefore
$\la(P_{j+1})_w,\delta^+\ra\leq0$, and $\delta^+w$ does not exist.

Finally we have to exclude the existence of $v^\vt w$. There are
two cases -- either $\la P_j^\vt ,w\ra<0$ or $\la P_j^\vt
,w\ra=0$. If $\la P_j^\vt ,w\ra<0$ we have $\la(P_{j+1})_w,v^\vt
\ra=\la P^\vt _j,v^\vt \ra=0$ by Lemma \ref{colv}, and we are
done. If $\la P_j^\vt ,w\ra=0$ and $v^\vt w$ existed, then its
image in $P_j^-$ under the $90^\circ$-rotation would be $vw$ -- a
contradiction because we have assumed that $vw$ does not exist.
\end{proof}

\section{Functorial properties}\label{BIFUN}

The essential data that determines the stable elementary group
$\EE(P)$ for a balanced polytope $P$ are
\begin{itemize}
\item[(i)] the matrix $\CB(P)=(\la F,v\ra)$ whose rows are
indexed by the column vectors $v$ of $P$ and whose columns are
indexed by the \emph{base facets} $F$ of column vectors of $P$,
and
\item[(ii)] the partial product structure on the set of
column vectors of $P$.
\end{itemize}

The following proposition shows that we are justified in saying
that polytopes $P$ and $Q$ are \emph{$E$-equivalent} if they
coincide in these data.

\begin{proposition}\label{EssData}
Let $P$ and $Q$ be balanced polytopes, and $R$ a ring. Suppose
there exists a map $\mu:\Col(P)\to\Col(Q)$ such that the following
conditions hold for all $v,w\in\Col(P)$:
$$
\emph{(i)}\quad\la P_w,v\ra=\la
Q_{\mu(w)},\mu(v)\ra\qquad\text{and}\qquad\emph{(ii)}\quad
\mu(vw)=\mu(v)\mu(w)\text{ if $vw$ exists.}
$$
\begin{itemize}
\item[(a)]
Then the assignment $x_v^\lambda\mapsto x_{\mu(v)}^\lambda$
induces a homomorphism
$$
\St(R,\mu):\St(R,P)\to\St(R,Q).
$$
\item[(b)]
Suppose that for all $v\in\Col(P)$ and all
$w\in\Col(Q)\setminus\mu(\Col(P))$ the following holds:
$-\mu(v)\neq w$, and none of the products $\mu(v)w$ and $w\mu(v)$
exists. (This condition obviously holds if $\mu$ is surjective.)
Then one has induced homomorphisms
$\EE(R,\mu):\EE(R,P)\to\EE(R,Q)$, and $K_2(R,\mu):K_2(R,P)\to
K_2(R,Q)$. $\EE(R,\mu)$ is surjective if $\mu$ is so.
\item[(c)]
If $\mu$ is bijective, then
$$
\St(R,P)\approx\St(R,Q),\quad \EE(R,P)\approx\EE(R,Q),\quad
K_2(R,P)\approx K_2(R,Q).
$$
\end{itemize}
\end{proposition}

\begin{proof}
(a) The mapping $\mu$ extends  to doubling spectra of $P$ and $Q$
and induces a map $\mu:\Col(\Pp)\to\Col(\Qq)$ as follows.

If $P$ is doubled with respect to $P_v$, then we also double $Q$
along $Q_{\mu(v)}$ and extend $\mu$ by setting
$\mu(v^-)=\mu(v)^-$, $\mu(v^\vt)=\mu(v)^\vt$,
$\mu(\delta^+)=\delta^+$, $\mu(\delta^-)=\delta^-$.

If $Q$ is doubled with respect to a facet $G$ that is not of type
$Q_\mu(v)$ for some $v\in\Col(P)$, then we put $\mu(v)=\mu(v)^-$.

Using Remark \ref{newcol}(b) one checks easily that the extension
of $\mu$ again satisfies the conditions (i) and (ii).

Condition (i) implies that the product $vw$ exists if and only if
$\mu(v)\mu(w)$ exists. Furthermore it follows from Proposition
\ref{maincrit}(d) that $v=-w$ if and only if $\mu(v)=-\mu(w)$. In
conjunction with (ii) this guarantees the compatibility of the
assignment $x_v^\lambda\mapsto x_{\mu(v)}^\lambda$ with the
Steinberg relations.

(b) We must check that $\St(R,\mu)$ maps the center of $\St(R,P)$
to the center of $\St(R,Q)$. This certainly holds if $x_v^\lambda$
for $v\in\Col(\Pp)$ commutes with every $x_w^\lambda$ for all
$w\in \Col(\Qq)\setminus\mu(\Col(\Pp))$. This is somewhat tedious,
but, in view of Remark \ref{newcol}(b), straightforward to check.
Hence the claim on induced homomorphisms follows from Theorem
\ref{ucext}, and that on surjectivity is then trivial.

(c) is obvious in view of Proposition \ref{centker}.
\end{proof}

\begin{remark}\label{EDRem}
(a) If the matrix $\CB(P)$ has pairwise different rows, then it
determines the partial product structure completely, since we can
identify the product $vw$ from $\CB(P)$ if it exists. Therefore,
every other polytope $Q$ such that $\CB(P)=\CB(Q)$ (up to a
suitable bijection $\mu:\Col(P)\to\Col(Q)$) has the same stable
Steinberg group. We will use this in Section \ref{POLYG}.

Since we evaluate column vectors only against base facets (and
\emph{not} against \emph{all} facets) to form $\CB(P)$, it is not
clear whether the product $vw$ can always be identified from
$\CB(P)$.

(b) Though we cannot prove $K_2$-functoriality for all maps $\mu$
discussed in Proposition \ref{EssData}, it is useful to note the
$\St$-functoriality, since it implies $K_i$-functoriality for
$i\ge 3$; see \cite{BrG5}.

(c) Proposition \ref{EssData} allows one to study polyhedral
$K$-theory as a functor also in the polytopal argument. The map
$\mu$ should be considered as a $K$-theoretic morphism from $P$ to
$Q$.
\end{remark}

As an easy application of Proposition \ref{EssData} one obtains
$$
K_2(R,P\times Q)=K_2(R,P)\oplus K_2(R,Q).
$$
for a ring $R$ and every pair of balanced polytopes $P$ and $Q$.
The analogous equations hold for $\St$ and $\EE$. In fact, one
observes that $\Col(P\times Q)$ is a disjoint union of $\Col(P)$
and $\Col(Q)$ (where $P$ and $Q$ are considered as lattice
subpolytopes of $P\times Q$), and that the column vectors coming
from $P$ are parallel to the base facets of those coming from $Q$,
and vice versa. These properties extend to doubling spectra.

Finally we want to point out that the $K$-theoretic groups only
depend on the projective toric variety associated with a polytope
$P$. The {\it normal fan} ${\cal N}(P)$  of a finite convex (not
necessarily lattice) polytope $P\subset\RR^n$ is defined as the
complete fan in the dual space $(\RR^n)^*=\Hom(\RR^n,\RR)$ given
by the system of cones
$$
\bigl(\{\phi\in(\RR^n)^*\mid\max_P(\phi)=F\},\ F\ \text{a face
of}\ P\bigr).
$$
Two polytopes $P,Q\subset\RR^n$ are called {\it projectively
equivalent} if ${\cal N}(P)={\cal N}(Q)$. In other words, $P$ and
$Q$ are projectively equivalent if and only if they have the same
dimension, the same combinatorial type, and the faces of $P$ are
parallel translates of the corresponding ones of $Q$.

\begin{remark} Let $P$ and $Q$ be {\it very ample} polytopes in the
sense of \cite[\S5]{BrG1}. This means that for any vertex $v\in P$
the affine semigroup $-v+(C_v\cap\ZZ^n)\subset\ZZ^n$ is generated
by $-v+\L_P$, where $C_v$ is the cone in $\RR^n$ spanned by $P$ at
$v$ (we assume that $P\subset\RR^n$ and $\gp(S_P)=\ZZ^{n+1}$), and
similarly for $Q$. Then ${\mathcal N}(P)={\mathcal N}(Q)$ if and
only if the projective toric varieties $\Proj(k[P])$ and
$\Proj(k[Q])$ are naturally isomorphic for some field $k$. These
varieties are normal, but not necessarily projectively normal
\cite[Example 5.5]{BrG1}
\end{remark}

Projectively equivalent polytopes $P$ and $Q$ have the same set of
column vectors: $\Col(P)=\Col(Q)$ (see \cite{BrG1}), and the
identity map on this set satisfies the condition of Proposition
\ref{EssData}.

\begin{proposition}\label{ProjEquiv}
If $P$ and $Q$ are projectively equivalent balanced polytopes then
$\St(R,P)\approx\St(R,Q)$ and $\EE(R,P)\approx\EE(R,Q)$.
\end{proposition}

\begin{remark}\label{coliso}
Proposition \ref{ProjEquiv} does not generalize to unstable groups
of elementary automorphisms. In fact, it is shown in \cite[Remark
5.3] {BrG1} that
$\Gamma_{\CC}(\Delta_1)\not=\Gamma_{\CC}(2\Delta_1)$, the
difference occurring exactly between $\EE_{\CC}(\Delta_1)$ and
$\EE_{\CC}(2\Delta_1)$.
\end{remark}

\section{Stable groups for balanced polygons}\label{POLYG}

First we classify all balanced polygons up to $E$-equivalence (and
in some cases up to projective equivalence), introducing
representatives for each class. Then we use this classification to
describe the corresponding stable groups of elementary
automorphisms in very explicit terms.

\begin{theorem}\label{dim2}
For a balanced polygon $P$ there are only the following
possibilities (up to affine-integral equivalence), each of which
can in fact be realized and constitutes a $E$-equivalence class:
\begin{itemize}
\item[(a)]
$P$ is a multiple of the unimodular triangle $P_a=\Delta_2$. Hence
$\Col(P)= \{\pm u,\pm v,\allowbreak\pm w\}$ and the column vectors
are subject to the obvious relations,
\item[(b)]
$P$ is projectively equivalent to the trapezoid
$P_b=\conv\bigl((0,0),(0,2),(1,1),\allowbreak(0,1)\bigr)$, hence
$\Col(P)=\{u,\pm v,w\}$ and the relations in $\Col(P)$ are $uv=w$
and $w(-v)=u$,
\item[(c)]
$\Col(P)=\{u,v,w\}$ and $uv=w$ is the only relation,
\item[(d)]
$\Col(P)$ has any prescribed number of column vectors, they all have the
same base edge (clearly, there are no relations between them),
\item[(e)]
$P$ is projectively equivalent to the unit lattice square $P_e$,
hence $\Col(P)=\{\pm u,\pm v\}$ with no relations between the
column vectors,
\item[(f)]
$\Col(P)=\{u,v\}$ so that $P_u\neq P_v$ with no relations in
$\Col(P)$.
\end{itemize}
\end{theorem}

\begin{proof}
It follows immediately from Proposition \ref{EssData} that the
polytopes in each class are $E$-equivalent, and that polytopes
from  different classes are not $E$-equivalent. Now let $P$ be a
balanced polygon.
\medskip

\noindent {\it Case 1}. There are column vectors $u,v\in\Col(P)$
such that $w=uv$ exists.

We claim that in this case $\Col(P)\subset\{\pm u,\pm v,\pm w\}$.
There exists a facet $F$ of $P$ such that $\la F,v\ra=\delta>0$.
Thus we have the following table whose entries are given by the
``scalar products'' of the vectors $u,v,w$ with the facets
$P_u,P_v,F$:
$$
\begin{pmatrix}
-1&1&\gamma\\
0&-1&\delta\\
-1&0&\gamma+\delta
\end{pmatrix}
$$
Note that $u,v,w$ have non-negative heights with respect to all
other facets.

Evidently $u$ and $v$ form a basis of the lattice $\ZZ^2$, and
every column vector $c$ is a linear combination of $u$ and $v$,
$c=\alpha u+\beta v$. Since $\la P_u,c\ra,\la
P_v,c\ra\in\{0,\pm1\}$, and at most one of these numbers can be
negative, one easily checks that $\alpha\in\{-1,0,1\}$ and the
following implications hold:
$$
\alpha=0\implies\beta\in\{\pm1\},\quad
\alpha=1\implies\beta\in\{0,1\}, \quad
\alpha=-1\implies\beta\in\{0,-1,-2\}.
$$
But $\alpha=-1$, $\beta=-2$ is impossible, since $\la F,c\ra\le-2$
in this case. We have now shown that indeed $\Col(P)\subset\{\pm
u,\pm v,\pm w\}$.

If $u,v,w$ are the only column vectors, then we are in case (c)
and can stop the discussion, since this case is possible: see the
proof of Lemma 8.3, where we have considered the polygon
$$
P_c=\conv\bigl((0,0),(3,0),(1,2),(0,1)\bigr).
$$
So suppose at least one of $-u,-v,-w$ is a further column vector.

Assume that $-u\in\Col(P)$. Then clearly $\gamma=0$, and $F$ is
parallel to $u$. After a change of coordinates we can assume
$u=(0,-1)$, $v=(1,0)$, so that $w=(1,-1)$. Since $P_u$ is parallel
to $v$, a parallel translation moves $P$ into a position where the
line segment $\conv\bigl((0,0),(1,0)\bigr)$ is contained in
$\partial P$ and $(1,0)$ is a vertex. The point $(0,1)$ must also
lie in $P_v$, which is parallel to $w$. Since $F$ is parallel to
$u$, the lower end point of $F$ must belong to $P_u$, and we can
sketch part of $P$ as in Figure \ref{BalPol}.
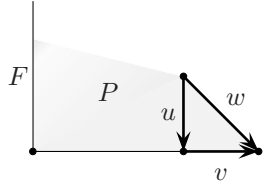
\begin{figure}[htb]
\begin{center}
 \psset{unit=1cm}
 \def\vertex{\pscircle[fillstyle=solid,fillcolor=black]{0.05}}
\begin{pspicture}(-2.5,-0.0)(1.5,2.0)
 \footnotesize
 \pspolygon[fillstyle=gradient,gradangle=-135, gradbegin=light,%
       gradend=white,linecolor=white](-2,0)(-2,1.5)(0,1)
 \pspolygon[fillstyle=solid, fillcolor=verylight,%
      linecolor=verylight](-2,0)(0,1)(1,0)
 \psline(-2,0)(0,0)
 \psline(-2,2.0)(-2,0)
 \rput(-2,0){\vertex}
 \rput(1,0){\vertex}
 \rput(0,1){\vertex}
 \rput(0,0){\vertex}
 \psline[style=fatline]{->}(0,1)(1,0)
 \psline[style=fatline]{->}(0,0)(1,0)
 \psline[style=fatline]{->}(0,1)(0,0)
 \rput(-1,0.8){$P$}
 \rput(-2.2,1.0){$F$}
 \rput(0.5,-0.3){$v$}
 \rput(-0.2,0.5){$u$}
 \rput(0.7,0.7){$w$}
\end{pspicture}
\end{center}
\caption{Part of the polygon $P$} \label{BalPol}
\end{figure}
The ``upper'' end point of $F$ must lie in $P_{-u}=P_v$ and so $P$
has exactly the three facets $P_u,P_v,F$. It follows that $P$ is a
multiple of the lattice unit triangle, and therefore
$\Col(P)=\{\pm u,\pm v, \pm w\}$. This is case (a).

If $-w$ is a column vector, then $-(\gamma+\delta)=-1$, and
consequently $\gamma=0$, $\delta=1$. Moreover, $-w$ and, hence,
$u,v,w$ have height $0$ with respect to all the other facets. Then
$-u$ is a column vector, and we are again in case (a).

Therefore, in addition to $u,v,w$ only $-v$ can be a column vector
(unless we are in case (a)). Then $F$ is the base facet of $-v$,
and this implies $\delta=1$, $\gamma=0$ (since $P$ is balanced),
so that $F$ is again parallel to $u$. In this case $P$ is
projectively equivalent to the polytope $P_b$, since all the
facets except $P_v$ and $P_{-v}=F$ must be parallel to $v$. (See
Figure \ref{ProdCol} where we have sketched an affine-integrally equivalent
polytope.)

\noindent{\it Case 2.} None of the products $uv$ ($u,v\in\Col(P)$)
exists.

There are two possibilities: either all the column vectors share
the same base edge or there are $u,v\in\Col(P)$ with $P_u\neq
P_v$. In the first situation it is clear that $P$ can have an
arbitrary number of column vectors: just consider the balanced
quadrangles
$$
P_{d,t}=\conv\bigl((-1,0),(3t,0),(2t,1),(0,2)\bigr),\ \ t\in\NN.
$$
In this situation $\Col(P_{d,t})=\{(s,-1)\ |\ s\in[0,t]\}$.
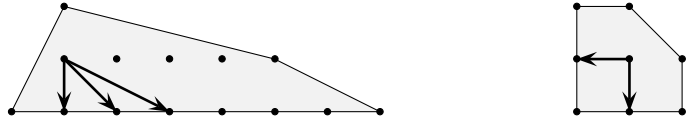
\begin{figure}[htb]
\begin{center}
 \psset{unit=0.7cm}
 \def\vertex{\pscircle[fillstyle=solid,fillcolor=black]{0.07}}
\begin{pspicture}(-1,-0.0)(6,2.0)
 \footnotesize
 \pspolygon[style=fyp](-1,0)(6,0)(4,1)(0,2)
 \multirput(-1,0)(1,0){8}{\vertex}
 \multirput(0,1)(1,0){5}{\vertex}
 \rput(0,2){\vertex}
 \psline[style=fatline]{->}(0,1)(1,0)
 \psline[style=fatline]{->}(0,1)(0,0)
 \psline[style=fatline]{->}(0,1)(2,0)
\end{pspicture}
\qquad\qquad\qquad
\begin{pspicture}(0,0)(2,2)
 \footnotesize
 \pspolygon[style=fyp](0,0)(2,0)(2,1)(1,2)(0,2)
 \multirput(0,0)(1,0){3}{\vertex}
 \multirput(0,1)(1,0){3}{\vertex}
 \multirput(0,2)(1,0){2}{\vertex}
 \psline[style=fatline]{->}(1,1)(0,1)
 \psline[style=fatline]{->}(1,1)(1,0)
\end{pspicture}
\end{center}
\caption{The polytopes $P_{d,t}$ and $P_f$} \label{ThreeColSame}
\end{figure}

Now assume $P_u\neq P_v$ for some $u,v\in\Col(P)$.

By Lemma \ref{colv} and Proposition \ref{maincrit}(a) we have $\la
P_v,u\ra=\la P_u,v\ra=0$. We claim that there are only two
possibilities: either $\Col(P)=\{u,v\}$ or $\Col(P)=\{\pm u,\pm
v\}$.

First we show that any element of $\Col(P)$ is parallel to either
$P_u$ or to $P_v$. In fact, if there is $w\in\Col(P)$ which is not
parallel to one of these edges then either $\la P_v,w\ra>0$ or
$\la P_u,w\ra>0$ (for otherwise $w$ would share the base edges
with both $u$ and $v$ which is impossible by Lemma \ref{colv}).
But then Proposition \ref{maincrit}(a) implies that one of the
products $wu$ or $wv$ exists -- a contradiction.

In particular we see that $\Col(P)\subset\{\pm u,\pm v\}$. One
only needs to show that if $-u\in\Col(P)$ then $-v\in\Col(P)$ as
well. In fact, if $-u\in\Col(P)$, then all the facets of $P$
different from $P_u$ and $P_{-u}$ are parallel to $u$, and $P_u$
and $P_{-u}$ are parallel to $v$. So $P$ is a parallelogram. Since
$u$ and $v$ span the lattice $\ZZ^2$, the parallelogram $P$ is
projectively equivalent to the unit square (up to affine-integral
equivalence). In this case $\Col(P)=\{\pm u,\pm v\}$.

That the case $\Col(P)=\{u,v\}$ is also possible is shown by the
balanced pentagon
$$
P_f=\conv\bigl((0,0),(2,0),(2,1),(1,2),(0,2)\bigr)
$$
whose only column vectors are $(-1,0)$ and $(0,-1)$.
\end{proof}

After the classification of balanced polygons up to
$E$-equivalence it remains to compute the stable elementary
groups. In the following theorem we use block matrix notation in a
self-explanatory way.

\begin{theorem}\label{sta2gp}
For a ring $R$ we have the group isomorphisms:
\begin{equation}
\tag{a}\EE(R,P_a)=\E(R),
\end{equation}
\begin{equation}
\tag{b}\EE(R,P_b)=
\begin{pmatrix}
\E(R)&\End_R(\oplus_{\NN}R)\\
&\\
0&\E(R)
\end{pmatrix},
\end{equation}
\begin{equation}
\tag{c}\EE(R,P_c)=
\begin{pmatrix}
\E(R)&\End_R(\oplus_{\NN}R)&
\Hom_R(\oplus_{\NN}R,R)\\
&\\
0&\E(R)&\Hom_R(\oplus_{\NN}R,R)\\
&\\
0&0&1
\end{pmatrix},
\end{equation}
\begin{equation}
\tag{d} \EE\bigl(R,P_{d,t}\bigr)=
\begin{pmatrix}
\E(R)&\Hom_R(\oplus_{\NN}R,R^t)\\
0&\text{\bf Id}_t
\end{pmatrix},\ \ t\in\NN,
\end{equation}
\begin{equation}
\tag{e}\EE(R,P_e)=\E(R)\times\E(R),
\end{equation}
\begin{equation}
\tag{f}\EE(R,P_f)=
\begin{pmatrix}
\E(R)&\Hom_R(\oplus_{\NN}R,R)\\
&\\
0&1
\end{pmatrix}
\times
\begin{pmatrix}
\E(R)&\Hom_R(\oplus_{\NN}R,R)\\
&\\
0&1
\end{pmatrix}.
\end{equation}
\end{theorem}

\begin{proof}
The first equation is just the definition of $\E(R)$.

It is enough to prove the second, third and fourth equations --
the last two follow from them since the unit square is the direct
product of two line segments, and since $P_f$ is $E$-equivalent to
$P_{d,1}\times P_{d,1}$.

Let $E_b$, $E_c$ and $E_{d,t}$ denote the groups on the right hand
side of the three equations. By $A=\{1',2',\ldots\}$ and
$B=\{1'',2'',\ldots\}$ we denote two mutually disjoint copies of
$\NN$, which are also disjoint from the original $\ZZ_+$, and
introduce the following index sets:
\begin{align*}
I_b&=\bigl(A\times A\bigr)\cup\bigl(B\times
A\bigr)\cup\bigl(B\times B\bigr)\\
I_c&=\bigl(A\times A\bigr)\cup\bigl(B\times A\bigr)\cup\bigl
(B\times B\bigr)\cup\bigl(\{0\}\times A\bigr)\cup\bigl(\{0\}\times
B\bigr),\\
I_{d,t}&=\bigl(A\times A\bigr)\cup\bigl(\{0,1,\ldots,t-1\}\times
A\bigr), \quad t\in\NN.
\end{align*}
By identifying $A$ and $B$ with free $R$-bases of two copies of
$\bigoplus_{\NN}R$, say $\FF_A$ and $\FF_B$, and $\{0,1,\ldots,t-1\}$
with that of $R^t$, we can view the
groups $E_b$, $E_c$, and $E_{d,t}$ as the corresponding subgroups
of $\Aut(\FF_A\oplus\FF_B)$, $\Aut_R(\FF_A\oplus \FF_B\oplus R)$
and $\Aut_R(\FF_A\oplus R^t)$ in an obvious way. In particular,
these groups consist of elementary automorphisms of free
$R$-modules with respect to their distinguished bases.

For each of the three cases we denote by $e_{ij}^\lambda$,
$\lambda\in R$ the associated elementary automorphism, $(i,j)$
belonging to the corresponding index set. They are subject to the
standard Steinberg relations.

Similar arguments as in the proof of Lemma \ref{specquad} (and
Proposition \ref{EssData}) show that $\St(R,P_b)$ and
$\St(R,P_{d,t})$ can be assumed to be generated by symbols
$x_{ij}^\lambda$ where $\lambda\in R$, $i\neq j$ and,
respectively, $(i,j)\in I_b$ or $(i,j)\in I_{d,t}$. (Lemma
\ref{specquad} itself is about $\St(R,P_c)$.) Moreover, these
symbols are subject to the same Steinberg relations as the
$e_{ij}^\lambda$. Therefore we have natural surjective group
homomorphisms
$$
\phi_b:\St(R,P_b)\to E_b,\quad\phi_c:\St(R,P_c)\to E_c,\quad
\phi_{d,t}:\St(R,P_{d,t})\to E_{d,t}.
$$
By Proposition \ref{centker} we only need to show the equations
$$
\Ker(\phi_b)=K_2(R,P_b),\quad\Ker(\phi_c)=K_2(R,P_c),\quad
\Ker(\phi_{d,t})=K_2(R,P_{d,t}).
$$
First we show that $Z(E_b)=Z(E_c)=Z(E_{d,t})=1$ ($Z$ is for the
center). Since $Z(\EE(R))=1$, the centers $Z(E_b)$, $Z(E_c)$ and
$Z(E_{d,t})$ are in the subgroups generated by $e_{ij}^\lambda$,
$\lambda\in R$, where the $(i,j)$ are respectively from $B\times
A$, $\bigl(B\times A\bigr)\cup\bigl(\{0\}\times
A\bigr)\cup\bigl(\{0\}\times B\bigr)$ and
$\{0,1,\ldots,t-1\}\times A$.

We have the group epimorphism $E_c\to E_b$ given by deleting the
last column and row. Therefore, if $Z(E_b)=1$ then $Z(E_c)$ is
actually concentrated in the smaller subgroup, generated by
$e_{ij}^\lambda$ with $(i,j)\in\bigl(\{0\}\times
A\bigr)\cup\bigl(\{0\}\times B\bigr)$.

Now we are done by the observation that for arbitrary
homomorphisms $B_0\in\End(\bigoplus_{\NN}R)$,
$B_1,B_2\in\Hom_R(\bigoplus_{\NN}R,R)$ and
$B\in\Hom_R(\bigoplus_{\NN}R,R^t)$ the equations
$$
\begin{pmatrix}
A_1&0\\
0&A_2
\end{pmatrix}\times
\begin{pmatrix}
1&B_0\\
0&1
\end{pmatrix}=
\begin{pmatrix}
1&B_0\\
0&1
\end{pmatrix}\times
\begin{pmatrix}
A_1&0\\
0&A_2
\end{pmatrix},
$$
$$
\begin{pmatrix}
A_1&0&0\\
0&A_2&0\\
0&0&1
\end{pmatrix}\times
\begin{pmatrix}
1&0&B_1\\
0&1&B_2\\
0&0&1
\end{pmatrix}=
\begin{pmatrix}
1&0&B_1\\
0&1&B_2\\
0&0&1
\end{pmatrix}\times
\begin{pmatrix}
A_1&0&0\\
0&A_2&0\\
0&0&1
\end{pmatrix}
$$
and
$$
\begin{pmatrix}
A&0\\
0&{\text{\bf Id}_t}
\end{pmatrix}\times
\begin{pmatrix}
1&B\\
0&{\text{\bf Id}_t}
\end{pmatrix}=
\begin{pmatrix}
1&B\\
0&{\text{\bf Id}_t}
\end{pmatrix}\times
\begin{pmatrix}
A&0\\
0&{\text{\bf Id}_t}
\end{pmatrix}
$$
simultaneously for all $A_1,A_2,A\in\EE(R)$ imply that $B_0$, $B_1$, $B_2$
and $B$ are actually zero homomorphisms.

It follows that $K_2(R,P_b)\subset\Ker(\phi_b)$,
$K_2(R,P_c)\subset\Ker(\phi_c)$ and $K_2(R,P_{d,t})\subset
\Ker(\phi_{d,t})$.

To derive the opposite inclusions we apply the same arguments as
in Proposition \ref{centker}. They go through once we show that
the analogue of the second part of the Claim in that proof remains
true when we change $\pi$ there by the unstable versions of each
of the homomorphisms $\phi_b$, $\phi_c$ and $\phi_{d,t}$.

More precisely, for each natural number $j$ we introduce the sets
$I_b^{(j)}$, $I_c^{(j)}$, $I_{d,t}^{(j)}$ which are finite
versions of $I_b$, $I_c$ and $I_{d,t}$, defined via the sets
$A^{(j)}=\{1',\ldots,j'\}$ and $B^{(j)}=\{1'',\ldots,j''\}$. They
give rise to the corresponding unstable subgroups
$E^{(j)}_b\subset E_b$, $E^{(j)}_c\subset E_c$,
$E^{(j)}_{d,t}\subset E_{d,t}$. We have surjective group
homomorphisms
$$
\phi_b^{(j)}:\St(R,P_b)_j\to E_b^{(j)},\quad\phi_c:\St(R,P_c)_j\to
E_c^{(j)},\quad\phi_{d,t}:\St(R,P_{d,t})_j\to E_{d,t}^{(j)},
$$
whose sources are the appropriate unstable Steinberg groups. One
should notice that for successive indices $j$ these Steinberg
groups may {\it not} correspond to successive members in the fixed
doubling spectra -- there may be big intervals of intermediate
members.

Assume $j>1$. We let
$$
{\mathfrak U}_b,{\mathfrak V}_b\subset\St(R,P_b)_j,\quad
{\mathfrak U}_c,{\mathfrak V}_c\subset\St(R,P_c)_j,\quad
{\mathfrak U}_{d,t},{\mathfrak V}_{d,t}\subset\St(R,P_{d,t})_j
$$
denote the subgroups which are the same for these unstable Steinberg groups
as ${\mathfrak U}^{i+1}$ and ${\mathfrak V}^{i+1}$ for $\St_R(P_{i+1})$ in the
proof of Proposition \ref{centker}.

What we want to show is that the restrictions
$$
\phi^{(j)}_b|_{{\mathfrak U}_b},\quad\phi^{(j)}_b|_{{\mathfrak V}_b},\quad
\phi^{(j)}_c|_{{\mathfrak U}_c},\quad\phi^{(j)}_c|_{{\mathfrak V}_c},\quad
\phi^{(j)}_{d,t}|_{{\mathfrak U}_{d,t}},\quad\phi^{(j)}_{d,t}|_{{\mathfrak
V}_{d,t}}
$$
are all injective group homomorphisms. We have the natural embeddings
$$
E_b^{(j)}\subset E_{2j}(R),\quad E_c^{(c)}\subset E_{2j+1}(R),\quad
E_{d,t}^{(j)}\subset E_{j+t}(R).
$$
Under these embeddings the generators of ${\mathfrak U}_b$
(${\mathfrak U}_c$, ${\mathfrak U}_{d,t}$) are sent to standard
elementary matrices $e_{st}^\lambda$ with the same $t$. Similarly,
the generators of ${\mathfrak V}_b$ (${\mathfrak V}_c$,
${\mathfrak V}_{d,t}$) are sent to standard elementary matrices
$e_{st}^\lambda$ with the same $s$. This, clearly, gives the
result.
\end{proof}

As the reader may guess Theorem \ref{sta2gp} implies that the
$K_2$ for balanced polygons is always either the usual $K_2$ or
twice $K_2$. We postpone the proof of this fact to \cite{BrG5}
where we treat all higher polyhedral $K$-groups simultaneously,
the polyhedral Milnor $K$-group being a basic ingredient in the
higher theory. Speaking loosely, all higher syzygies between
elementary automorphisms come from unit simplices -- a stable
higher version of Theorem \ref{plg} in the polygonal case.

However, such a nice matrix theoretical interpretation of the
stable groups of elementary automorphisms as in Theorem
\ref{sta2gp} is no longer possible for higher dimensional
polytopes as explained in the example below. This makes the
computation of polyhedral $K$-groups challenging.

\begin{example}\label{nomatrix}
In the proof of Corollary \ref{assoc} we have discussed the
pyramid $P$ over the unit square (see Figure \ref{ThreeCol}). This
polytope has $8$ column vectors and $5$ facets, each of which is
the base facet for at least one column vector. In the table $(\la
F,v\ra)$ the rows corresponding to the column vectors whose base
facet is the unit square have \emph{two} entries $1$. Moreover,
every edge represents a column vector. All of this makes it impossible
to construct a
matrix representation in which elementary automorphisms are mapped to
standard elementary matrices. (Though, we have the canonical embedding
$\EE_R(P)\to\E_5(R)$, $5=\#\L_P$.)
\end{example}

\end{document}